\documentclass[journal]{IEEEtran}
\IEEEoverridecommandlockouts
\usepackage{amsmath,amsthm,amssymb,graphicx,mathtools,hyperref,multicol,booktabs,float,mathrsfs,textcomp,xcolor,amsfonts,comment}
\usepackage{tikz,pgfplots}

\newtheorem{remark}{Remark}
\newtheorem{definition}{Definition}
\newtheorem{problem}{Problem}
\newtheorem{theorem}{Theorem}
\newtheorem{lemma}{Lemma}

\newtheorem{corollary}{Corollary}
\newtheorem{assumption}{Assumption}

\definecolor{orange2}{RGB}{255,165,0}
\definecolor{violet2}{rgb}{0.93, 0.51, 0.93}
\definecolor{green2}{rgb}{0.2, 0.75, 0.2}

\newcommand{\sto}{ \rm{s.t.} }

\newcommand{\norm}[1]{\left\lVert #1 \right\rVert}

\newcommand{\rank}[1]{\text{rank}\left(#1\right)}

\newcommand{\Lag}{\mathcal{L}}
\newcommand{\zstar}{z^{\star}}
\newcommand{\xstar}{x^{\star}}
\newcommand{\lstar}{\lambda^{\star}}
\newcommand{\xstart}{x^{\star\top}}
\newcommand{\lstart}{\lambda^{\star\top}}
\newcommand{\controller}{\mathcal{K}}
\newcommand{\plant}{\mathcal{P}}
\newcommand{\R}{\mathbb{R}}
\newcommand{\akaZ}{\mathcal{Q}}
\newcommand{\akaB}{\mathcal{R}}
\newcommand{\g}{g}
\newcommand{\argmin}[1]{\underset{#1}{\operatorname{arg}\,\operatorname{min}}\;}

\begin{document}

\title{A new framework for constrained optimization\\ via feedback control of Lagrange multipliers}
\author{V. Cerone, \and S. M. Fosson, \and S. Pirrera, \and D. Regruto
\thanks{The authors are with the Dipartimento di Automatica e Informatica, Politecnico di Torino,
    corso Duca degli Abruzzi 24, 10129 Torino, Italy;
    e-mail: vito.cerone@polito.it, sophie.fosson@polito.it, simone.pirrera@polito.it, diego.regruto@polito.it.
}}

\maketitle

\begin{abstract}
The continuous-time analysis of iterative algorithms for optimization has a long-standing history.
This work introduces a novel framework for equality-constrained optimization based on control theory.
The central concept is to design a feedback control system in which the Lagrange multipliers serve as the control inputs while the output represents the constraints. This system converges to a stationary point of the constrained optimization problem through suitable regulation. Concerning the Lagrange multipliers, we explore two control laws: proportional-integral control and feedback linearization. These choices lead to a variety of different methods. We rigorously develop the related algorithms, analyze their convergence theoretically, and present several numerical experiments that demonstrate their effectiveness compared to the state-of-the-art approaches.
\end{abstract}


\section{Introduction}\label{sec:IN}
First-order iterative algorithms are prevalent in convex and non-convex optimization and machine learning to handle large-scale datasets and leverage parallel processing architectures. Iterative algorithms for optimization are discrete-time (DT) dynamical systems that update the estimate of the optimization variables at each iteration. Their continuous-time (CT) counterparts, obtained by considering infinitesimal step sizes, are described by differential equations whose analysis can provide a deeper theoretical understanding, particularly concerning stability and convergence rate.

A paradigmatic example is the gradient flow, defined by the equation $\dot{x}=-\nabla f(x)$, where $f:\R^n\to\R$ is a differentiable, unconstrained cost function that we aim to minimize. Gradient flow is the CT version of gradient descent and proximal minimization algorithms achieved through forward and backward Euler discretization,  respectively; see, e.g., \cite[Sec. 4.1.1]{par13}. The study of gradient flow is particularly relevant in deep learning, where gradient descent methods have seen practical success in training despite lacking a solid theoretical understanding; see, e.g., \cite{sax14,elk21} and references therein. Additionally, works \cite{su16,mue19} propose CT analysis of Nesterov accelerations for gradient descent.

In this work, we consider the more challenging problem of constrained optimization. A significant amount of research focuses on CT methods based on Lagrange multipliers,  particularly primal-dual methods. For a general overview, we refer the reader to \cite[Chapter 15]{lue_book}. The main CT approach to constrained optimization is the primal-dual gradient dynamics (PDGD), which was introduced in \cite{kos56,arr58}.

In \cite{qu19}, the authors examine the exponential stability of PDGD when minimizing strongly convex, smooth cost functions with linear equality constraints. This analysis is extended to non-smooth composite optimization in \cite{din19,dhi19} by using a proximal augmented Lagrangian, and to non-convex stochastic optimization in \cite{che19}. Regarding convex composite optimization with linear equality constraints, the work \cite{fra18} illustrates a CT model for the alternating direction method of multipliers (ADMM \cite{boy10}).
Additionally, prior research proposes modifying the gradient flow to account for equality constraints. The key idea in this approach is to build the descent direction $\dot{x}$ as a combination of a projected gradient and a Gauss-Newton direction, which drives the solution toward the feasible set; see, e.g.,  \cite{yam80,sch00,zho07}.
Finally, works \cite{sch00,zho07,qu19}  also consider inequality constraints by extending first-order conditions to Karush-Kuhn-Tucker (KKT) conditions \cite{lue_book}.

This brief review emphasizes that much of the literature focuses on the CT analysis of existing algorithms rather than the development of new CT algorithms for optimization.
An exception is the recent work \cite{all24}, which develops CT dynamical systems that solve constrained optimization problems. The authors synthesize a feedback controller, exploiting tools from the theory of control barrier functions (CBF), that guarantees forward invariance and stability of the feasible set. These properties enforce safety when the feasible set represents the safe operation of a plant.

This work introduces a new CT framework for convex and non-convex constrained optimization, focusing on equality constraints. The proposed framework adopts an original feedback control approach: we start from the first-order necessary conditions for minima to build a CT dynamical system where the vector of Lagrange multipliers represents the control input. The output of this system corresponds to the constraints, which we regulate accordingly. Several control laws can be employed for the Lagrange multipliers to achieve the desired regulation, leading to a family of control-based first-order methods. 

This work has two main contributions. The first one is developing and analyzing a control-theoretic framework that enables the synthesis of new first-order optimization algorithms using feedback control design techniques. We call this framework controlled multipliers optimization (CMO).

The second contribution is the specialization of CMO to specific control laws. We focus mainly on proportional-integral control and feedback linearization methods, denoted as PI-CMO and FL-CMO, respectively. We develop the corresponding algorithms, we study the conditions for their convergence in both convex and non-convex settings and the convergence speed. Through this analysis, we highlight the advantages of the proposed algorithms compared to state-of-the-art CT techniques for constrained optimization, with particular attention to references \cite{qu19,all24}.
Finally, we present several numerical experiments demonstrating the effectiveness of PI-CMO and FL-CMO, even beyond their theoretically established convergence conditions.
We pay special attention to comparing our approaches to state-of-the-art optimization algorithms regarding solution accuracy and numerical complexity.

We organize the paper as follows. In Section \ref{sec:PS}, we formulate the problem, review the theory of Lagrange multipliers, and describe the proposed control-theoretic framework. The two succeeding sections specialize the framework into two distinct control strategies. Specifically, Section \ref{sec:PI} develops the PI control method, demonstrates its convergence, and analyzes the convergence rate for strongly convex problems with linear constraints. Section \ref{sec:FL} introduces the feedback linearization method and proves its convergence for both strongly convex and non-convex problems. Section \ref{sec:NR} presents several numerical experiments that illustrate the practical effectiveness of the proposed methods in various applications. Finally, Section \ref{sec:CON} concludes the paper.
\section{Problem statement and proposed framework}\label{sec:PS}
We consider the constrained optimization problem
\begin{equation} \begin{split}
    \min_{x \in \R^n}& \, f(x) \\ &\sto \\ &h(x) = 0
    \label{opt}
\end{split} \end{equation}
where $f: \R^n \mapsto \R$ and $h: \R^n \mapsto \R^m$ are differentiable, possibly non-convex functions.

The Lagrangian of problem \eqref{opt} is
\begin{equation}\label{eq:lagrangian}
\Lag(x,\lambda)=f(x)+\lambda^\top h(x)
\end{equation}
where $\lambda\in\R^m$ is the vector of Lagrange multipliers. We report the following well-known theorem for self-consistency; see, e.g., \cite[Sec 11.3]{lue_book} for a complete overview.
\begin{theorem}[First-order necessary conditions]
\label{lagr_th}
Let $\xstar \in \R^n$ be a local minimum of $f$ such that $h(\xstar)=0$. Assume that $\xstar$ is regular, i.e., $\nabla h_1(\xstar), \dots,\nabla h_m(\xstar)$ are linearly independent. Then, there exists a unique  $\lstar \in \R^m$ such that $(\xstar,\lstar)$ is a saddle point of $\Lag(x,\lambda)$, i.e,
\begin{equation} \label{FO}
   \nabla f(\xstar) + J_h(\xstar)^\top\lstar   = 0
\end{equation}
where $J_h(x)\in\R^{m,n}$ is the Jacobian matrix of $h$ evaluated in $x$.
\end{theorem}
In the rest of the paper, we call stationary point any couple $(\xstar, \lstar)$, with $\xstar\in\R^n$ and $\lstar \in\R^m$, which satisfies \eqref{FO} and the constraints, i.e.,
\begin{equation}\label{eq:stationary}
\begin{split}
  &\nabla f(\xstar) + J_h(\xstar)^\top\lstar   = 0\\
  &h(\xstar)=0.
\end{split}
 \end{equation}
In this work, we focus on first-order algorithms, i.e., on methods that achieve stationary points. In the non-convex case, finding a stationary point does not guarantee local optimality; however, the literature is rich in first-order methods in large-scale optimization and machine learning thanks to their low complexity when compared to, e.g., second-order methods, which require computing and storing Hessian matrices as well; see, e.g., \cite{ber99,lue_book,good16}.

\subsection{Proposed framework: feedback control of Lagrange multipliers}\label{sub:lagrange}
In this subsection, we illustrate the proposed framework to find a stationary point of problem \eqref{opt} by building a suitable CT dynamical system with controlled Lagrange multipliers.

\begin{figure}
    \centering
    \includegraphics[width=0.93\columnwidth]{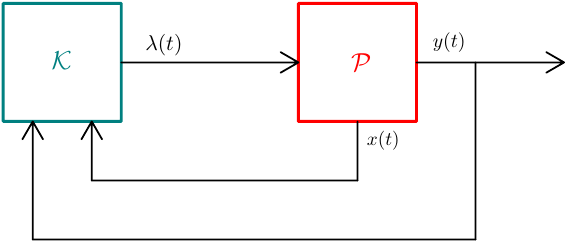}
    \caption{Structure of the proposed feedback control approach. The state and the output of $\plant$ defined in \eqref{uncontrolled} are fed back to the controller $\controller$, whose output is the vector of the Lagrange multipliers.}
    \label{fig:fcs}
\end{figure}

Let us define the dynamical system $\plant$ with state $x(t) \in \R^n$, input $\lambda(t) \in \R^m$ and output $y(t) \in \R^m$,  described by the equations
\begin{equation}
\label{uncontrolled}
    \plant: \left\{ \begin{split}
        \dot{x}(t) &= -\nabla f(x(t)) - J_h(x(t))^\top\lambda(t) \\
        y(t) &= h(x(t))
    \end{split} \right.
\end{equation}
The following result holds.
\begin{lemma}\label{lem1}
An equilibrium point $(\xstar,\lstar)$ of $\plant$  is a stationary point of problem \eqref{opt} if and only if $h(\xstar)=0$.
\end{lemma}
\begin{proof}
By definition, an equilibrium point $(\xstar,\lstar)$ of $\plant$ satisfies equation \eqref{FO}. If $h(\xstar)=0$, then \eqref{eq:stationary} holds and $(\xstar,\lstar)$ is a stationary point to problem \eqref{opt}. Conversely, any stationary point satisfies \eqref{eq:stationary}, thus \eqref{FO}, then it corresponds to an equilibrium point $(\xstar,\lstar)$ of $\plant$ with  $h(\xstar)=0$.
\end{proof}
Lemma \ref{lem1} suggests we can compute a stationary point of problem \eqref{opt} by designing a suitable input $\lambda(t)$ that drives $\plant$ to converge to an equilibrium point and regulates the output to zero. A standard way to approach this regulation problem is to design a suitable feedback controller $\controller$. In Fig. \ref{fig:fcs}, we depict a general scheme;  $y(t)$  is the feedback signal to the input of $\controller$, possibly together with the state of $\plant$. We underline that $\plant$ is a representation of an optimization algorithm; therefore, the state is known as well as the output, which is different from physical systems where the observation of the state may be critical.

The goal of this work is to tackle the following problem.
\begin{problem}\label{problem}
Design a feedback controller $\controller$ for $\plant$ such that
\begin{equation}
 \begin{split}
 & \lim_{t\to\infty}x(t)= \xstar,\\
 &\lim_{t\to\infty}y(t)=0.
 \end{split}
\end{equation}
\end{problem}

In this work, we consider two possible design techniques for $\controller$ to solve Problem \ref{problem}: PI control and feedback linearization.

\subsection{Related work}
To our knowledge, controlling the Lagrange multipliers in constrained optimization is novel. In the literature, the use of control methods to analyze or develop optimization algorithms is only partly explored.

Regarding first-order unconstrained optimization, the works \cite{les16} and \cite{hu17} present a control theory interpretation of known algorithms such as gradient descent, heavy-ball, and Nesterov's accelerated methods. In particular, they show that these algorithms correspond to  DT feedback systems, where the current gradient is the control input. By using integral quadratic constraints, they analyze their convergence. Unlike our work, this framework does not envisage constrained optimization.

As mentioned in the introduction, the work \cite{all24} proposes a CT algorithm based on CBF that solves equality and inequality constrained optimization. The proposed method can also be interpreted as a continuous approximation of a projected gradient flow \cite{xia00,cor08}, which is used to analyze the stability.

A different research line studies the solution of equations via control. In the DT framework, \cite{men22} and \cite{men23} solve linear algebraic equations using iterative learning control and observer-based controller design, respectively. Instead, in \cite{bha07}, the authors build a CT-controlled system whose output is regulated to solve $g(x)=0$, where $g:\R^n\mapsto \R^n$ is a vector function. By choosing appropriate control Lyapunov functions, they retrace standard iterative methods, such as Newton-Raphson and conjugate gradient methods, and develop new variants.
Unlike our work, this framework considers equality-constrained problems with no cost function to minimize. However, one may argue that finding a stationary point corresponds to solving the first-order equations $(\nabla f(x) + \sum_{i=i}^m \lambda_i\nabla h_i(x))^\top$ together with $h(x)=0$ in the variables $x$ and $\lambda$, i.e., a system of dimension $(n+m)\times (n+m)$. In other terms, one can apply the methods proposed in \cite{bha07} to find the zeros of the vector function $g(x,\lambda)=\left[ (\nabla f(x) + \sum_{i=i}^m \lambda_i\nabla h_i(x))^\top, h(x)^\top\right]^\top$. However, as summarized in \cite[Table 1]{bha07}, this approach gives rise to second-order methods and continuous Newton algorithms, which require the inversion of the Jacobian of $g(x,\lambda)$. Therefore, the numerical complexity is prohibitive for large-scale problems.

Beyond the control-theoretic approach, PDGD \cite{qu19} is a valuable CT approach to solve problem \eqref{opt}. In the next section, we compare our approach to PDGD.

\section{First method: PI-CMO}\label{sec:PI}
\begin{figure}
    \centering
    \includegraphics[width=0.93\columnwidth]{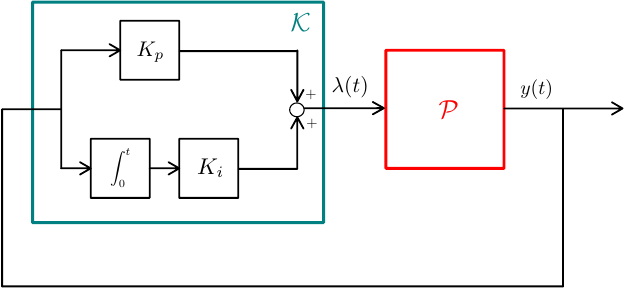}
    \caption{First method: PI-CMO. We feedback the output $y(t)$ to the controller $\controller$, which applies proportional an integral actions.}
    \label{fig:PI_scheme}
\end{figure}
A possible strategy to design $\controller$ for Problem \ref{problem} is to apply a PI action on $y(t) = h(x(t))$. The PI control is widespread in industrial and engineering applications, thanks to its effectiveness in regulating many processes by tuning two parameters.

In our framework, the PI control law is as follows:
\begin{equation}
    \lambda(t) = K_p y(t) + K_i \int_{0}^{t} y(\tau)\,\rm{d} \tau
    \label{input_pi}
\end{equation}
where $K_p\in\R$ and $K_i\in\R$ are the coefficients of the proportional and integral terms, respectively. We depict the corresponding feedback scheme in Fig. \ref{fig:PI_scheme}.

As a consequence, $\controller$ is a dynamical system described by the differential equation
\begin{equation} \begin{split}
    \dot{\lambda}(t) &=K_p   \frac{\rm{d}}{\rm{d}t}  y(t) + K_i y(t)\\
    &= K_p J_h(x(t)) \dot{x}(t) + K_i y(t).
\end{split}\end{equation}

In the following, we drop the variable $t$ in long formulas, namely $x=x(t)$, $\dot{x}=\dot{x}(t)$, $\lambda=\lambda(t)$ and $\dot{\lambda}=\dot{\lambda}(t)$.
Given the definition of $\dot{x}$ in \eqref{uncontrolled}, the closed-loop dynamics is
\begin{equation} \begin{split}
     \dot{x} &= -\nabla f(x) - J_h(x)^\top\lambda\\
    \dot{\lambda} &= -K_p J_h(x) \left[\nabla f(x) + J_h(x)^\top\lambda\right]+ K_i h(x).
    \label{cl}
\end{split} \end{equation}
or, equivalently,
\begin{equation} \begin{split}
     \dot{x} &= -\nabla_x\Lag(x,\lambda)\\
    \dot{\lambda}&= -K_p J_h(x) \nabla_x\Lag(x,\lambda)+ K_i \nabla_{\lambda}\Lag(x,\lambda).
    \label{cl2}
\end{split} \end{equation}
In the following, we refer to the dynamics \eqref{cl} as to PI-CMO dynamics. For any equilibrium point $(\xstar,\lstar)$, from \eqref{cl} we get the condition $h(\xstar)=0$. Therefore, according to Lemma \ref{lem1}, each equilibrium point corresponds to a stationary point in the case of PI control.

\begin{remark}
In equation \eqref{input_pi}, we consider scalar $K_p$ and $K_i$ to reduce the number of design parameters and simplify the convergence analysis reported in section \ref{sub:PI_conv}. Extension to the most general case where $K_p$ and $K_i$ are matrix gains in $\R^{m, m}$ is straightforward.
\end{remark}
\begin{remark} 
It is worth noticing that by considering a purely integral control, i.e., $K_p=0$, we obtain the dynamics
\begin{equation} \begin{split}
      \dot{x} &= -\nabla_x\Lag(x,\lambda)\\
    \dot{\lambda}&=  K_i \nabla_{\lambda}\Lag(x,\lambda),
    \label{PDGD}
\end{split} \end{equation}
which corresponds to PDGD defined in \cite[Eq. (2a)-(2b)]{qu19}, where the symbol $\eta$ is used instead of $K_i$. In other words, we can interpret PDGD as applying an integral control to the Lagrange multipliers, which recast PDGD into the proposed control-theoretic framework. In contrast, the proposed PI-CMO extends PDGD to a novel family of CT optimization algorithms.
\end{remark}

\subsection{Global exponential convergence of PI-CMO for strongly convex problems }\label{sub:PI_conv}
This section proves that the PI-CMO dynamics converges globally and exponentially in a convex setting with affine constraints, using an appropriate Lyapunov function. Additionally, we analyze the convergence rate. We adopt the same assumptions as those in \cite{qu19} to enable a thorough comparison with PDGD.
\begin{assumption}\label{ass:sc}
 $f$ is strongly convex and twice differentiable.
\end{assumption}
\begin{assumption}\label{ass:aff}
 $h$ is affine, i.e., $h(x)=Cx+d$, $C\in\R^{m,n}$, $d\in\R^m$.
Moreover, $C$ is full rank and there exist $0<\alpha_1\leq \alpha_2$ such that
\begin{equation}\label{CCposdef}
 \alpha_1 I \preceq CC^\top \preceq \alpha_2 I.
\end{equation}
\end{assumption}
The assumption that $C$ is full rank guarantees that $Cx+d=0$ has solutions and that there are no linearly dependent constraints. Additionally, under this assumption, any point $x \in \R^n$ is regular. Therefore, all assumptions required by Theorem~\ref{lagr_th} are satisfied.
We define
\begin{equation}
 z(t):=(x(t)^\top,\lambda(t)^\top)^\top
\end{equation}
and
\begin{equation}
 \zstar:=(\xstart,\lstart)^\top
\end{equation}
is the equilibrium point of \eqref{cl}, which corresponds to a saddle point of $\Lag(x,\lambda)$. Assumptions \ref{ass:sc} and \ref{ass:aff} guarantee the existence and the uniqueness of $\zstar$.

According to \cite[Lemma 1]{qu19}, there exists a symmetric $B(x)\in\R^{n,n}$ satisfying $\beta_1 I \preceq B(x) \preceq \beta_2 I$ for some $0< \beta_1<\beta_2$ such that
\begin{equation}\label{lem:qu}
 \nabla f(x)-\nabla f(\xstar) = B(x)(x-\xstar).
\end{equation}
The constants $\beta_1$ and $\beta_2$ satisfy
\begin{equation}\label{eq:beta12_hessian}
    \beta_1 I\leq \nabla^2 f(\xstar + t(x-\xstar))  \leq \beta_2 I, \quad \forall t \in [0,1].
\end{equation}
where $\nabla^2 f(x)$ is the Hessian matrix of $f(x)$.
Equation \eqref{eq:beta12_hessian} relates the constants $\beta_1$ and $\beta_2$ with the second-order information of $f(x)$. If $f(x)$ is quadratic and strongly convex, namely $f(x)=\frac{1}{2}x^\top W x$ with $W\succ 0$, then $\beta_1$ and $\beta_2$ are the minimum and maximum eigenvalues of $W$, respectively. 
If $f(x)$ is not a quadratic function, we can evaluate $\beta_1$ and $\beta_2$  as outlined in \cite{Monn08}, which presents a method for estimating tight bounds on the eigenvalues of $\nabla^2 f(x)$.

In the following, we write $B=B(x)$ for brevity.

\begin{theorem}[Global exponential convergence of PI-CMO]
\label{theo:pi}
Let assumptions \ref{ass:sc} and \ref{ass:aff} hold. Given $K_i>0$ and $K_p>0$, if
\begin{equation}\label{kp}K_p< \frac{2K_i}{\beta_2},\end{equation}
 then there exist real positive constants $c_1$ and $c_2$ such that
\begin{equation}
 \|x(t)-\xstar\|_2\leq c_1e^{-\frac{1}{2}\mu t},~~~\|\lambda(t)-\lstar\|_2\leq c_2e^{-\frac{1}{2}\mu t}
\end{equation}
where
\begin{equation}\label{tau}
\mu =\min\left\{ K_p \alpha_1,2\beta_1 -\frac{K_p}{K_i}\beta_1\beta_2\right\}>0.
\end{equation}
\end{theorem}

\begin{proof}
First of all, we define the candidate Lyapunov function
\begin{equation}\label{def:lyap}
V\big(z(t)\big)=\big(z(t)-\zstar\big)^\top P \big(z(t)-\zstar\big)
\end{equation}
where
\begin{equation}\label{def:P}
P:=\left(\begin{array}{cc}
                   K_i I_n&0\\
                    0&I_m
                   \end{array}\right)\in\R^{m+n,m+n}.
            \end{equation}
If we prove that
\begin{equation}\label{lyap}
 \dot{V}\big(z(t)\big)\leq -\mu V\big(z(t)\big)
\end{equation}
for some $\mu>0$, then the theorem statement holds. Therefore, in the following, we focus on conditions that guarantee \eqref{lyap}.

We start with some preliminary computations.
Since $\nabla_x\Lag(\xstar,\lstar)=0$, $\nabla_{\lambda}\Lag(\xstar,\lstar)=0$ and $J_h(x)=C$, we have
\begin{equation}\label{gradx}
 \begin{split}
  \nabla_x\Lag(x,\lambda)&=\nabla_x\Lag(x,\lambda)-\nabla_x\Lag(\xstar,\lstar)\\
  &=\nabla f(x)-\nabla f(\xstar)+J_h(x)^{\top}\lambda-J_h(\xstar)^{\top}\lstar\\
  &=\nabla f(x)-\nabla f(\xstar)+C^{\top}(\lambda-\lstar)
 \end{split}
\end{equation}
and
\begin{equation}\label{gradl}
 \begin{split}
  \nabla_{\lambda}\Lag(x,\lambda)&=\nabla_{\lambda}\Lag(x,\lambda)-\nabla_{\lambda}\Lag(\xstar,\lstar)\\
  &=Cx+d-(C\xstar+d) = C(x-\xstar).
 \end{split}
\end{equation}
By using \eqref{cl2}, \eqref{lem:qu}, \eqref{gradx} and \eqref{gradl}, we obtain
\begin{equation}
 \begin{split}
\dot{z}&(t)=\left(\begin{array}{c}
                    \dot{x}(t)\\
                    \dot{\lambda}(t)\\
                   \end{array}\right)\\
                   &=\left(\begin{array}{c}
                    -B(x-\xstar)-C^\top(\lambda-\lstar)\\
                    K_iC(x-\xstar) -K_pC[B(x-\xstar)+C^\top(\lambda-\lstar)]
                   \end{array}\right)\\
                  &= \left(\begin{array}{cc}
                    -B&-C^\top\\
                    K_i C -K_pCB&-K_p CC^\top
                   \end{array}\right)
                  \big(z(t)-\zstar\big).
 \end{split}
\end{equation}
Let us define \begin{equation}
G:=\left(\begin{array}{cc}
                    -B&-C^\top\\
                    K_i C -K_pCB&-K_p CC^\top
                   \end{array}\right)
            \end{equation}
so that $$\dot{z}(t)=G \big(z(t)-\zstar\big).$$
Then,
\begin{equation}
\begin{split}
 \dot{V}\big(z(t)\big)&=\dot{z}(t)^{\top}P\big(z(t)-\zstar\big)+\big(z(t)-\zstar\big)^\top P\dot{z}(t)\\
 &=\big(z(t)-\zstar\big)^\top\big(G^\top P+P G\big)\big(z(t)-\zstar\big).
\end{split}
\end{equation}
Therefore, a sufficient condition for  $\dot{V}\big(z(t)\big)\leq -\mu V\big(z(t)\big)$, see \eqref{lyap},
is
\begin{equation}\label{cond}
-G^{\top} P -P G-\mu P \succeq 0.
\end{equation}
As a result, our next goal is to establish sufficient conditions for \eqref{cond}.
We compute
\begin{equation}\label{pg}
PG=\left(\begin{array}{cc}
                    -K_i B&-K_i C^\top\\
                    K_i C -K_pCB&-K_p CC^\top
                   \end{array}\right)
\end{equation}
while $G^{\top} P=(PG)^\top$. Hence,
\begin{equation}\label{gg}
\begin{split}
 &-G^{\top} P -P G-\mu P\\
 &= \left(\begin{array}{cc}
                    2K_i B-K_i\mu I& K_pB C^\top\\
                    K_p C  B& 2 K_p CC^\top-\mu I
                   \end{array}\right)\\
%                     &= \left(\begin{array}{cc}
%                     2\rho B-\rho \mu I&\left[(K_i -\rho) I -  K_pB\right]C^\top\\
%                     C\left[(K_i-\rho) I - K_pB\right]& K_p CC^\top-\mu I
%                    \end{array}\right)\\
 &\succeq \left(\begin{array}{cc}
                    2K_i B-K_i \mu I& K_pB C^\top\\
                    K_pC B & K_p CC^\top
                   \end{array}\right)
                   \end{split}
\end{equation}
where the last step derives from $2K_p CC^\top-\mu 
 I\succeq  CC^\top$ which holds for $\mu\leq  K_p \alpha_1$.
Since $CC^\top \succ 0$ is invertible from \eqref{CCposdef} in Assumption \ref{ass:aff}, we can apply the Schur complement argument to conclude that the matrix
\begin{equation*} \left(\begin{array}{cc}
                    2K_i B-K_i \mu I&K_pBC^\top\\
                    K_p CB& K_p CC^\top
                   \end{array}\right)
\end{equation*}
is positive semidefinite if and only if
\begin{equation}\label{schur}
2K_i B-K_i \mu I-K_p BC^\top \frac{1}{K_p}(CC^\top)^{-1}K_p C B\succeq 0.
\end{equation}
Moreover, since $CC^\top$ is invertible, then $C^{\top}(CC^\top)^{-1}C\preceq I$. Therefore, a sufficient condition for \eqref{schur} is
\begin{equation}\label{sufficient}
(2K_i I- K_p B)B\succeq K_i \mu I.
\end{equation}
Under assumption \eqref{kp}, 
\eqref{sufficient} holds if
\begin{equation}
 (2K_i-K_p\beta_2) \beta_1\geq K_i \mu
\end{equation}
which is equivalent to
\begin{equation}
 \mu \leq 2\beta_1-\frac{K_p}{K_i}\beta_1\beta_2.
\end{equation}
%From \eqref{rho}, this bound is non-negative.
This completes the proof.
\end{proof}

\begin{remark}
According to \eqref{def:lyap}-\eqref{def:P}, the considered Lyapunov function is
 $$V(x,\lambda)=K_i\|x-\xstar\|_2^2+\|\lambda-\lstar\|_2^2.$$
We highlight the role of the parameter $K_i$, which tunes the weight assigned to the terms in $x$ relative to the terms in $\lambda$. According to Theorem \ref{theo:pi}, a larger $K_i$ increases the bound on $\mu$, meaning that the weight given to the terms in $x$ must be sufficiently large compared to the terms in $\lambda$. We remark that in the convergence proof of PDGD in \cite{qu19}, the Lyapunov function is defined on a non-diagonal matrix $P$, which complicates the interpretation of the parameters.
\end{remark}

According to \cite{qu19}, PDGD is globally exponentially convergent with rate $\frac{1}{2}\mu_{PDGD}$
\begin{equation}\label{mupdgd}
\mu_{PDGD}=\min\left\{\frac{\eta\alpha_1}{4\beta_2},\frac{\alpha_1\beta_1}{4\alpha_2}\right\}%\leq \frac{\alpha_1\beta_1}{4\alpha_2}.
\end{equation}
Since $\eta$ is a design parameter, we can set $\eta\geq\frac{\beta_1 \beta_2}{\alpha_2}$ it so that $\mu_{PDGD}$ saturates to $\frac{\alpha_1\beta_1}{4\alpha_2}$.
The following result holds.
\begin{corollary}\label{cor:PI}
For $\epsilon\in \left(0,1-\frac{\alpha_1}{8\alpha_2}\right)$, let \begin{equation}
K_p = \epsilon\frac{2K_i}{\beta_2}~\text{ and }~
K_i> \frac{1}{\epsilon}\frac{\beta_1\beta_2}{8\alpha_2}.
\end{equation}
Then
\begin{equation}\label{taua}
\mu >\mu_{PDGD}
\end{equation}
i.e., PI-CMO enjoys a faster convergence rate with respect to PDGD \cite{qu19}.
\end{corollary}
\begin{proof}
By replacing $K_p = \epsilon\frac{2K_i}{\beta_2}$ in \eqref{tau}, we obtain
$$\mu=\min\left\{2\epsilon K_i\frac{\alpha_1}{\beta_2},2\beta_1(1-\epsilon)\right\}.$$
Now, $\epsilon\in \left(0,1-\frac{\alpha_1}{8\alpha_2}\right)$ implies $2\beta_1(1-\epsilon)>\frac{\alpha_1\beta_1}{4\alpha_2}$, while $K_i> \frac{1}{\epsilon}\frac{\beta_1\beta_2}{8\alpha_2}$ implies $2\epsilon K_i\frac{\alpha_1}{\beta_2}>\frac{\alpha_1\beta_1}{4\alpha_2}$.
This proves the statement because $\mu_{PDGD}\leq\frac{\alpha_1\beta_1}{4\alpha_2}$.
\end{proof}
In Sec. \ref{sub:first}, we validate this result about the enhanced convergence speed through numerical simulations.
\subsection{Illustrative example}\label{sub:kikp}
To complete the analysis, we present an example that illustrates how the tuning of $K_i$ and $K_p$ may affect the convergence rate. We also compare the results to the case $K_p=0$, representing PDGD.
Let us consider the univariate optimization problem
\begin{equation}
 \min_{x\in\R}\frac{1}{2} w x^2~~ \text{s.t.}~~~ x=0
\end{equation}
where $w>0$.
The PI-CMO dynamics for this problem is
\begin{equation} \begin{split}
     \dot{x} &= -wx - \lambda\\
    \dot{\lambda} &= (K_i- K_pw)x-K_p \lambda.
    \label{cl_ex}
\end{split} \end{equation}
This is a second-order CT linear time-invariant system
\begin{equation}
\begin{pmatrix}\dot{x}\\ \dot{\lambda}\end{pmatrix}=\mathrm{A}
\begin{pmatrix}x\\ \lambda \end{pmatrix}
\end{equation}
with
\begin{equation}
 \mathrm{A}:=\left(\begin{array}{cc}
 -w & -1 \\
 K_i- K_pw & -K_p
 \end{array}\right).
\end{equation}
The eigenvalues of $\mathrm{A}$ are
\begin{equation}\label{eigA}
\frac{-(K_p+w)\pm \sqrt{(K_p+w)^2-4K_i}}{2}.
\end{equation}
For any $K_i>0$, the eigenvalues are either real and negative or complex with negative real parts. In particular, if $K_i\geq (K_p+w)^2/4$, $K_i$ contributes only to the imaginary part, therefore it does not impact on the convergence rate.

According to  \cite{qu19}, although $\eta>0$ can be arbitrarily large for PDGD, increasing $\eta$ beyond a certain threshold does not lead to a faster decaying rate.
In our example, if we choose $K_p=0$, we obtain PDGD with $\eta=K_i$ and, from \eqref{eigA}, we see that if $K_i\geq w^2/4$, then $K_i$ has no impact on the real parts of the eigenvalues. This explains the observation in \cite{qu19}.

On the other hand, in PI-CMO, we can tune $K_p>0$ to enhance the convergence rate, provided that the conditions of Theorem \ref{theo:pi} are satisfied, which is a benefit with respect to PDGD.

\subsection{PI-CMO in non-convex quadratic optimization with linear constraints}\label{sub:quad}
To conclude this section, we analyze some properties of PI-CMO for quadratic optimization with linear constraints. In particular, we prove that optimization problems exist with non-convex cost functions in which PI-CMO converges to a stationary point while PDGD is divergent.

We consider
\begin{equation}
\begin{split}
&\min_{x\in\R^n} \frac{1}{2}x^\top W x\\
&\text{s.t.}\\
&Cx+d=0
\end{split}
\end{equation}
where $W\in \R^{n,n}$ is symmetric and $C\in\R^{m,n}$ and $d\in\R^m$, with invertible $CC^\top$.

Since the cost function is quadratic and the constraints are linear, PI-CMO corresponds to the linear time-invariant system
\begin{equation}\label{forquadratic}
\begin{pmatrix}
                    \dot{x}\\
                    \dot{\lambda}\\
                   \end{pmatrix}
                   = \mathrm{A}
                  \begin{pmatrix}
                    x\\
                    \lambda\\
                   \end{pmatrix}+\begin{pmatrix}
                    0\\
                    K_id\\
                   \end{pmatrix}
\end{equation}
with
\begin{equation}
    \mathrm{A}:=\begin{pmatrix}
       -W&-C^\top\\
                    K_i C -K_pCW&-K_p CC^\top 
    \end{pmatrix}
\end{equation}

Therefore, if $\mathrm{A}$ is Hurwitz, then the system is asymptotically stable, and by construction, the output $y(t)=Cx(t)+d$ is regulated to zero.
Thus, for PI-CMO $W$ does not need to be positive definite, i.e., the system does not need to be strongly convex.
As an example, let us consider $W=\text{diag}(1,-1)$ and $C=(0,2)$. Since $W$ is indefinite, the quadratic cost function is not convex, while the constrained problem has a unique minimum.
The eigenvalues of the corresponding dynamic matrix are $-1$ and $\frac{1-4K_p\pm\sqrt{(1-4K_p)^2-16K_i}}{2}$. If $K_p=0$, all the eigenvalues have a positive real part for all $K_i$. In other terms, PDGD is always unstable. Instead, for $K_p>\frac{1}{4}$, all the eigenvalues have negative real part for all $K_i>0$.

In conclusion, there exist problems with non-convex functions where PI-CMO converges to the minimum as long as we provide a suitable tuning of $K_p$. In contrast, PDGD is divergent for any hyperparameter choice. This observation encourages future study of PI-CMO in non-convex optimization.

\section{Second method: FL-CMO}\label{sec:FL}
In this section, we resort to feedback linearization as detailed, e.g., in \cite{isi95}, to design the controller $\controller$ introduced in Section \ref{sec:PS}. Moreover, we study the conditions under which the controlled dynamics is stable and the algorithm converges to the desired solution.

We organize the section as follows. Sec. \ref{sub:fl_basics} summarizes the key concepts and results necessary for the development of the proposed algorithm and its stability analysis. In particular, we overview the non-interacting control problem and its solution. In Sec. \ref{sub:boh}, we apply the non-interacting control framework to Lagrange multipliers, which leads to the definition of the FL-CMO algorithm. The last subsections are devoted to the convergence analysis.

\subsection{Feedback linearization basics}\label{sub:fl_basics}
In this subsection, we review some basic concepts of the feedback linearization theory to design the control of $\lambda$ and to analyze the convergence of the resulting dynamical system.
\begin{definition}\emph{(Lie derivative, \cite[Sec. 1.2]{isi95}, \cite[Sec. 13.2]{kha02})}.

    Let $F: \mathbb{R}^n \mapsto \mathbb{R}^n$ be a vector field and $H_i: \mathbb{R}^n \mapsto \mathbb{R}$. The Lie derivative of $H_i$ along $F$ is
    \begin{equation}\label{lie_der}
        L_{F} H_i (x)= \nabla H_i(x)^\top F(x)\in\R.
    \end{equation}
    By defining $L_F^0 H_i(x) = H_i(x)$, for $k=1,2,\dots$, we have
    \begin{equation}
        L_F^k H_i(x) = (\nabla  L_F^{k-1}H_i(x))  ^\top F(x).
    \end{equation}
\end{definition}

As illustrated, e.g., in \cite{isi95,kha02}, feedback  linearization can be applied to input-affine nonlinear dynamical systems of the form
\begin{equation} \label{sys_nl} \begin{split}
    \dot{x} &= F(x) +  G(x) u\\
    y &= H(x)
\end{split} \end{equation}
where $x(t) \in \R^n$, $u(t) \in \R^m$, and $y(t) \in \R^m$, $F:\mathbb{R}^n  \mapsto \mathbb{R}^n$, $G: \mathbb{R}^n  \mapsto \mathbb{R}^{n,m}$, $H: \mathbb{R}^n  \mapsto \mathbb{R}^m$.
For a dynamical system of this kind, we recall the definition of the relative degree.

Let $G_j(x)$ be the $j$-th column of $G(x)$.
\begin{definition}\emph{(Relative degree, \cite[Sec. 5.1]{isi95}, \cite[Sec. 13.2]{kha02})}.

System \eqref{sys_nl} has a vector relative degree $r = (r_1,\dots,r_m)^\top$ at $\bar{x}\in\R^n$ if

    (a) for each $1 \leq i,j \leq m$, for each $k \in \mathbb{N}$ such that  $k< r_i-1$, and for all $x$ in a neighbourhood of $\bar{x}$

    \begin{equation} \label{rel_deg}
        L_{G_j} L_F^k H_i(x) = 0
    \end{equation}
    and

    (b) the matrix $\left[ L_{G_j} L_F^{r_i-1} H_i(x) \right]_{1 \leq i,j \leq m}\in\R^{m,m}$ is nonsingular at $x=\bar{x}$, i.e.,
    \begin{equation} \label{rel_deg2}
        \text{rank} \left( \left[ L_{G_j} L_F^{r_i-1} H_i(\bar{x}) \right]_{1 \leq i,j \leq m} \right) = m.
    \end{equation}
\end{definition}
For a single-input, single-output system, if the relative degree is $r$, the output $y$ and its derivatives up to $(r-1)$-th order do not depend on the
input, while the $r$-th order derivative does. For linear systems, the concept is equivalent to the relative degree of a transfer function, i.e., the difference between the degree of the denominator and the degree of the numerator.

Now, we recall the concept of zero dynamics.  \begin{definition}\emph{(Zero dynamics, \cite{byr84}, \cite[Sec. 4.3]{isi95})}.
Consider system \eqref{sys_nl}. Let $(\xi^\top,\eta^\top)^\top = \Phi(x)$ be a smooth change of coordinates such that $\xi = (y_1,\dot y_1,\dots y_1^{(r_1-1)},\dots,y_m,\dots,y_m^{(r_m-1)}) \in \R^{r}$ and $\eta \in \R^{n-r}$. Using this transformation, we obtain the normal form representation of \eqref{sys_nl}, i.e.,
\begin{equation}\begin{aligned}
    \dot \xi &= \phi(\xi, \eta) + \gamma(\xi,\eta) u \\
    \dot \eta &= q(\gamma, \eta),
\end{aligned}  \end{equation}
where $\phi: \R^{r}\times \R^{n-r} \rightarrow \R^r$, $\gamma: \R^{r}\times \R^{n-r} \rightarrow \R^{r,m}$, and $q: \R^{r}\times \R^{n-r} \rightarrow \R^{n-r}$. The zero dynamics of \eqref{sys_nl} is defined as the dynamics of the system $\dot \eta = q(0,\eta)$.
\end{definition}
The zero dynamics of a system $\mathcal{S}$ describes the internal behavior of $\mathcal{S}$ when the output is constrained to be identically zero through a suitable choice of input and initial conditions. We refer the reader to \cite{isi13} for further insights.

Finally, we state the non-interacting control problem and present its solution. 
\begin{definition}\emph{(Non-interacting control problem, \cite[Sec. 5.3]{isi95})}
    Let us consider system \eqref{sys_nl}. We say that a controller of the form
    \begin{equation}\label{uabx}
        u = \alpha(x) + \beta(x) v
    \end{equation}
    with $\alpha(x)\in\R^m$, $\beta(x)\in\R^{m,m}$ and $v\in\R^m$
    solves the non-interacting control problem if in the closed-loop system
        \begin{equation}\label{controlled_sys_nl} \left\{ \begin{split}
            \dot{x} &= F(x) + G(x) \alpha(x) + G(x) \beta(x) v \\ y &= H(x)
        \end{split} \right. \end{equation}
    each input $v_i$ affects only output $H_i$, for each $i=1,\dots,m$.
\end{definition}

The solution to the non-interacting control problem is established in \cite{isi95}; for completeness, we report this result in the following theorem.

 \begin{theorem}\emph{(Non-interacting control solution, \cite[Sec. 5.3]{isi95})}\label{th:non_interact}

        (a) The non-interacting control problem admits a solution if and only if system \eqref{sys_nl} has some vector relative degree $r$. In that case, the solution is given by
        \begin{equation}
            \alpha(x)= -A^{-1}(x) b(x), \qquad \beta(x) = A^{-1}(x)
        \end{equation}
        where
        \begin{equation}\label{Ab}
        \begin{split}
            &A(x) = \left[ L_{G_j} L_F^{r_i-1} H_i(x) \right]_{1 \leq i,j \leq m}\in\R^{m,m}\\
            &b(x) = \left[ L_F^{r_i} H_i(x) \right]_{i=1,\dots,m}\in\R^m.
        \end{split}
        \end{equation}
        Moreover, for each $i=1,2,\dots,m$, the input-output behavior between $v_i(t)$ and $y_i(t)$ is linear and described in the $s-$domain by the transfer function
        \begin{equation} \label{th_non_interact_linearized}
            \frac{1}{s^{r_i}} .
        \end{equation}
(b) If $v(t)$ stabilizes the system described by \eqref{th_non_interact_linearized} and the zero dynamics of \eqref{sys_nl} is
asymptotically stable, then the feedback control system is asymptotically stable.
    \end{theorem}

Part (a) of Theorem \ref{th:non_interact} provides an input-output global linearization of system \eqref{sys_nl}, which is decoupled in each component $i=1,\dots,m$. From part (b) of Theorem \ref{th:non_interact}, the stability of the controlled system also depends on the zero dynamics, which accounts for the non-observable part of the system.

\subsection{Non-interacting control of Lagrange multipliers}\label{sub:boh}
In this section, we apply the theory reported in Sec. \ref{sub:fl_basics} to design a feedback linearized control for the plant $\plant$ as defined in \eqref{uncontrolled}.

The key point is that $\plant$ has the input-affine structure of \eqref{sys_nl}, with
\begin{equation}\label{pezzi}
\begin{split}
F(x)=-\nabla f(x),~G(x)=-J_h(x)^\top,~ H(x)= h(x).
 \end{split}
\end{equation}
Therefore, we can apply Theorem \ref{th:non_interact} to obtain a decoupled controller if $\plant$ has some vector relative degree.
Let us analyze this point using the following assumption.

\begin{assumption} \label{ass:fl2}
    Every $x\in\R^n$ is a regular point for problem \eqref{opt}, i.e.,
    \begin{equation}
        \rank{J_h(x)} = m, \qquad \forall x \in \R^n.
    \end{equation}
\end{assumption}
In other terms, the Jacobian of the constraints is full row rank; this implies $m \leq n$, i.e., there are at least as many optimization variables as constraints. This assumption is quite standard in constrained optimization; see, e.g., \cite{nocedal2006numerical,sch00,all24}. In practice, it is met in many applications, including a broad class of system identification problems, such as in the example reported in Sec. \ref{sub:third}, optimal control problems \cite{dor95} and distributed optimization over networks \cite{hon18}.

In the case of affine constraints $h(x) = C x + d = 0$, Assumption \ref{ass:fl2} guarantees that at least one feasible solution exists. Conversely, in the case of non-convex constraints, Assumption \ref{ass:fl2} implies that the feasible set $\Omega=\{x\in\R^n: h(x)=0\}$ is an $(n-m)$-dimensional smooth manifold, i.e., locally similar to a $\R^{n-m}$ at all points. Consequently, $\Omega$ is not empty, and we can define a dynamical system whose trajectories evolve onto it, i.e., the zero dynamics.
We finally remark that Assumption \ref{ass:fl2} implies the regularity condition of Theorem \ref{lagr_th}.

\begin{lemma}\label{lemma_reldeg1}
Under Assumption~\ref{ass:fl2}, the system $\plant$ in \eqref{uncontrolled} has a vector relative degree $r =(1,1,\dots,1)^\top \in \R^m$.
\end{lemma}
\begin{proof}
Let us consider \eqref{sys_nl} with \eqref{pezzi}. 
It is sufficient to check \eqref{rel_deg2} with $r_i=1$ for each $i=1,\dots,m$ to prove the thesis. From \eqref{pezzi}, $H_i(x) = h_i(x)$ and $G_j(x) =
-\nabla h_j (x)$. Let us set $r_i=1$ in \eqref{rel_deg2}. Then,
        \begin{equation}
            L_{G_j} L_F^0 H_i(x) = L_{G_j} H_i(x) = \left( \nabla h_i(x) \right)^\top \left( -\nabla h_j(x)\right).
        \end{equation}
        Thus,
        \begin{equation}\begin{split}
&\text{rank} \left( \left[ L_{G_j} L_F^{0} H_i(x) \right]_{1\leq i,j\leq m} \right) = \\ \quad &= \text{rank} \left( \left[\left( \nabla h_i(x) \right)^\top \left( -\nabla h_j(x)\right)
\right]_{ij} \right) \\ \quad &= \text{rank} \left( - J_h(x)  J_h(x) ^\top \right)= m.
        \end{split} \end{equation}
        Hence, \eqref{rel_deg2} is satisfied for any $x$.
    \end{proof}
According to Lemma \ref{lemma_reldeg1}, part (a) of Theorem \ref{th:non_interact} holds for system \eqref{controlled_sys_nl} with \eqref{pezzi}.
Therefore, we resort to the non-interacting control solution given in Theorem \ref{th:non_interact}. Specifically, we define the static feedback control law
    \begin{equation} \label{decoupling_control}
        \lambda(t) = A(x)^{-1}(-b(x) + v(t))
    \end{equation}
    where $A(x)$ and $b(x)$ are given by \eqref{Ab} with $r_i=1$, i.e.,
     \begin{equation}\label{Abr1}
        \begin{split}
            &A(x) = \left[ L_{G_j} H_i(x) \right]_{1 \leq i,j \leq m}=- J_h(x)  J_h(x)^\top\in\R^{m,m}\\
            &b(x) = \left[ L_F H_i(x) \right]_{i=1,\dots,m}=- J_h(x) \nabla f (x) \in \R^{m}.
        \end{split}
        \end{equation}

    By applying the control law \eqref{decoupling_control}-\eqref{Abr1}, the relationship between the new input $v(t) \in \R^m$ and the output $y(t)=h(x(t)) \in \R^m$   is
    \begin{equation}\label{linearized_sys} y_i(t) = \int_{0}^{t} v_i(\mu)\, {\rm{d}} \mu,~~~i=1,\dots,m, \end{equation}
    which is linear and decoupled in each component $i$.

\begin{figure}
    \centering
    \includegraphics[width=0.93\columnwidth]{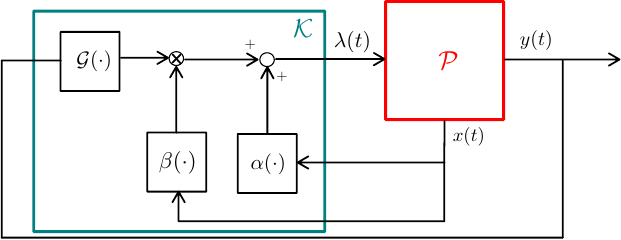}
    \caption{Second method: feedback linearization. We feedback the output $y(t)$ to the controller $\controller$ to compute $v(t)$ according to \eqref{linear_control}, and the state $x(t)$ to compute $u(t)$ according to \eqref{uabx}.}
    \label{fig:FL_scheme}
\end{figure}

    Next, we have to design $v(t)$ to regulate $y(t)$ to zero. A possible solution is to design the controller $\mathcal{G}$ in Fig. \ref{fig:FL_scheme}
\begin{equation}\label{linear_control}
v(t) = \mathcal{G}( y(t) )
\end{equation}
such that the closed-loop dynamics is asymptotically stable and $y(t) \to 0$.

In the following, we refer to the closed-loop system dynamics defined by equations \eqref{uncontrolled}, \eqref{decoupling_control} and \eqref{linear_control}, i.e., 
\begin{equation}\label{eq:flcmo}
    \dot{x}(t) = -\nabla f(x) + J_h^\top (J_h J_h^T)^{-1}(J_h\nabla f(x) + \mathcal{G}(y)),
\end{equation}
as to FL-CMO dynamics.

Given the single integral structure of \eqref{linearized_sys}, the simplest way to design   $\mathcal{G}$ is to consider $m$ static linear feedback controllers 
\begin{equation}\label{linear_control_details}
   v_i(t) = -K_i y_i(t) 
\end{equation} with $K_i>0$ for $i=1,\dots,m$. In fact, this leads to
\begin{equation}
    {y_i(t)} = \alpha_i e^{-K_i t}, \quad i=1,\dots,m
\end{equation}
where $\alpha_i$ is a constant that depends on the initial conditions. 
Several more possibilities are available for the design of $\mathcal{G}$, each leading to a different CT optimization algorithm in the FL-CMO family characterized by different properties. In Sec. \ref{sub:fifth}, we present a numerical example comparing different choices for $\mathcal{G}$ and draw additional considerations on this aspect.

\begin{remark}\label{rem:bottleneck}
    FL-CMO requires the inversion of $A(x)$, which is a computational bottleneck for large $m$. To address this point, we notice that $-A(x)=J_h(x) J_h(x)^\top$ is positive definite by construction. Consequently, we can efficiently address the inversion through Cholesky factorization. Alternatively, we can compute the QR factorization $J_h(x)^\top= QR$. This choice yields the Cholesky factorization of $-A$ while avoiding the computation of the product $J_h(x)J_h(x)^\top$; in fact, $-A = R^\top Q^\top Q R = R^\top R$.  The sparsity level of $J_h(x)$ determines which of the two methods is more effective.
\end{remark}

\begin{remark}
The FL-CMO dynamics \eqref{eq:flcmo} limited to the case of static $\mathcal{G}$ defined in  \eqref{linear_control_details} substantially corresponds to the systems studied in \cite{yam80,tan80}. The recent work \cite{all24} retrieves those systems and reframes them in the context of CBF-based feedback control, with an extension to the case of inequality constraints.

Interestingly, the derivation of the algorithm in \cite{all24} is entirely different from the one presented in this work. The intuition behind \cite{all24} is to add a drift to the standard gradient, aiming to keep the feasibility of the state via a control action. In particular, CBF theory is used to synthesize the controller. In the case of equality constraints, the admissible control set is a singleton that corresponds to \eqref{linear_control_details}; see \cite[Remark 4.2]{all24} for details.

The proposed FL-CMO extends the CBF approach to a family of controllers $\mathcal{G}$ beyond \eqref{linear_control_details}. Although a comprehensive study of the properties of different $\mathcal{G}$'s is beyond the purpose of this work, in Sec. \ref{sub:fifth}, we propose a numerical example that compares static and dynamic $\mathcal{G}$'s in a robust optimization problem to highlight the benefits of dynamic controllers.
\end{remark}

From the development of FL-CMO, we can observe that its performance depends on two factors. On the one hand, the rate of convergence to the feasible set $\Omega$ is completely determined by the user when designing the external controller $\mathcal{G}$. In other terms, the convergence to $\Omega$ is arbitrarily fast.
On the other hand, the zero dynamics determines the convergence speed to the optimal solution. Since the convergence to $\Omega$ is arbitrarily fast, the overall time required for convergence is determined in practice by the uncontrollable zero dynamics. Nevertheless, numerical experiments demonstrate that if the convergence rate to $\Omega$ is taken too fast, the differential equations defining the closed-loop system become stiff, rendering the numerical integration more challenging.

In Fig. \ref{fig:FL_scheme}, we summarize the feedback control scheme obtained via feedback linearization.

\subsection{Local convergence of FL-CMO}
In this section, we analyze the convergence of FL-CMO for possibly non-convex problems.

First, we notice that the FL-CMO dynamics converges to the feasible set $\Omega$ by construction. To prove the pointwise convergence, let us introduce the second-order sufficient conditions.

\begin{definition}\label{def:2nd}\emph{(Second-order sufficient conditions)}.
For problem \eqref{opt}, let $H_{xx}\Lag(x,\lambda)$ be the Hessian matrix of $\Lag(x,\lambda)$ with respect to $x$. We say that the second-order sufficient conditions hold at $(\xstar,\lstar)$ if \begin{equation}
v^\top H_{xx}\Lag(\xstar,\lstar) v >0
\end{equation}
for any $v\in\R^n$ such that $J_h(\xstar)v=0$.
\end{definition}
If the second-order sufficient conditions hold, $\xstar$ is a strict local minimum of \eqref{opt}.
We prove the following result of local convergence.
    \begin{theorem}[Local convergence of FL-CMO to local minima]\label{th:fl_convergence}
    Let the second-order sufficient conditions hold at $(\xstar,\lstar)$. Then the local minimum $\xstar$ is a locally asymptotically stable equilibrium point of the FL-CMO  dynamics.
    \end{theorem}

\begin{proof}
	Let us consider part (b) of Theorem \ref{th:non_interact}. First, we notice that $v(t)$ stabilizes \eqref{th_non_interact_linearized} by construction and regulates the output to zero, see \eqref{linear_control}.

	Therefore, it is sufficient to prove that the zero dynamics of \eqref{sys_nl} is asymptotically stable to obtain the asymptotic stability of the closed-loop dynamics \eqref{controlled_sys_nl}-\eqref{pezzi}-\eqref{linear_control}. We prove this fact in a neighbourhood of an equilibrium point $(\xstar, \lstar)$ of $\plant$, see \ref{sub:lagrange}.

To analyze the zero dynamics, we begin by defining the mapping $\Phi: \R^n\mapsto \R^n$ as
\begin{equation}\label{def:phi}
	\Phi(x) = \begin{pmatrix} h(x) \\ J_h^\perp(\xstar)(x- \xstar)   \end{pmatrix},
\end{equation}
where we define $J_h^\perp(x) \in \R^{n-m,n}$ as follows: its rows are an orthonormal basis for the null space of the rows of $J_h$. As a consequence, $J_h^\perp(x)J_h^\top(x)=0$ for all $x\in\R^n$.
The Jacobian matrix of $\Phi(x)$ is
\begin{equation}
	J_{\Phi}(x)= \begin{pmatrix}
	           J_h(x) \\ J_h^\perp(\xstar)
	          \end{pmatrix}.
\end{equation}
Since $J_h(x)$ has rank $m$ by Assumption \ref{ass:fl2} and $J_h^\perp(x)$ has orthogonal rows by definition, then $\rank{J_{\Phi}(x)} = n.$ Therefore, $\Phi$ is invertible and provides a suitable change of coordinates in the state space of $\plant$. In particular, we can express the transformed state $z=\Phi(x)$ as
\begin{equation}\label{def:z}
z=\begin{pmatrix} \xi \\ \eta \end{pmatrix}
\end{equation}
where $\xi=y\in\R^m$ corresponds to the output and $\eta\in\R^{n-m}$ represents the state of the zero dynamics of the system. We refer the reader to \cite[Sec 5.1]{isi95}.

Now, by exploiting the change of variables via $\Phi$, we write the system normal form and analyze the zero dynamics.

From \eqref{uncontrolled}, we have
$$\dot{x}=\g(x,\lambda):=-\nabla f(x) - J_h(x)^\top\lambda.$$
Then, the normal form is
\begin{equation} \label{normalform_1} \begin{split}
	\dot{z} &= \frac{\text{d}}{\text{d} t} \Phi(x)  =J_{\Phi}(x) \g(x,\lambda)=J_{\Phi}(\Phi^{-1}(z)) \g(\Phi^{-1}(z),\lambda).
\end{split} \end{equation}
We notice  that $\Phi(\xstar)=0$. Then, let us represent \eqref{normalform_1} through its Taylor expansion around $(\Phi(\xstar), \lstar)=(0,\lstar)$, i.e.,
\begin{equation} \label{expansion_normal_form}\begin{split}
	\dot{z} =  \akaZ z + \akaB (\lambda-\lstar) + o(z)+o(\lambda-\lstar).
\end{split} \end{equation}
Specifically,
\begin{equation} \begin{split}
	\akaZ &= \frac{\partial}{\partial z} \left[J_{\Phi}(\Phi^{-1}(z))\g(\Phi^{-1}(z),\lstar)\right]_{\vert{z=0}} \\
&= \frac{\partial}{\partial x} \left[J_{\Phi}(x)\g(x,\lstar)\right]_{\vert{x=\xstar}}  J_{\Phi^{-1}}(0)\\
\end{split} \end{equation}
Since $\g(\xstar,\lstar)=0$,
\begin{equation}
\begin{split}
&\frac{\partial}{\partial x} \left[J_{\Phi}(x)\g(x,\lstar)\right]_{\vert{x=\xstar}}=\\
&=\frac{\partial}{\partial x} \left[J_{\Phi}(x)\right]_{\vert{x=\xstar}}\g(\xstar,\lstar)+ J_{\Phi}(\xstar)\frac{\partial}{\partial x} \left[\g(x,\lstar)\right]_{\vert{x=\xstar}}\\
&= J_{\Phi}(\xstar)\frac{\partial}{\partial x} \left[\g(x,\lstar)\right]_{\vert{x=\xstar}}.
\end{split}
\end{equation}
Thus,
\begin{equation} \begin{split}
	\akaZ =-J_{\Phi}(\xstar)H_{xx}\Lag(\xstar,\lstar)J_{\Phi^{-1}}(0),
\end{split} \end{equation}
where
\begin{equation}\begin{split}
	H_{xx}\Lag(\xstar,\lstar) &=  -\frac{\partial}{\partial x} \left[\g(x,\lstar)\right]_{\vert{x=\xstar}}.
\end{split} \end{equation}
Moreover, by the inverse function theorem
\begin{equation}\begin{split}
	 J_{\Phi^{-1}}(0) &=  J_{\Phi}^{-1}(\xstar) = \begin{pmatrix} J_h^{\dagger}(\xstar),  & J_h^{\perp \top}(\xstar) \end{pmatrix}
\end{split}\end{equation}
where $J_h^{\dagger}(\xstar)=J_h^\top(\xstar)[J_h(\xstar) J_h^\top(\xstar)]^{-1}$.
In conclusion,
\begin{equation}
	\akaZ = -J_{\Phi}(\xstar) H_{xx}\Lag(\xstar,\lstar) \begin{pmatrix} J_h^{\dagger}(\xstar), & J_h^{\perp \top}(\xstar) \end{pmatrix}.
\end{equation}
On the other hand, given $\frac{\partial}{\partial \lambda} \g(\xstar,\lambda)_{\vert{\lambda=\lstar}} = -J_h^\top(\xstar)$, we have
\begin{equation} \begin{split}
	\akaB &=  \frac{\partial}{\partial \lambda}  \left( J_{\Phi}(\xstar)\g(\xstar,\lambda))\right\vert_{\lambda=\lstar}     \\
&= -J_{\Phi}(\xstar)J_h^\top(\xstar) =- \begin{pmatrix} J_h(\xstar) J_h^\top(\xstar)  \\ J_h^\perp(\xstar) J_h^\top(\xstar)  \end{pmatrix} = \\
&=  -\begin{pmatrix} J_h(\xstar) J_h^\top(\xstar)  \\ 0 \end{pmatrix}.
\end{split} \end{equation}
Next, we obtain the zero dynamics by setting $\xi=0$ and considering the last $n-m$ equations of \eqref{normalform_1}:
\begin{equation}  \begin{split}
\dot \eta &=   -J_h^\perp(\xstar) H_{xx}\Lag(\xstar,\lstar) \begin{pmatrix} J_h^{\dagger}(\xstar), & J_h^{\perp \top}(\xstar) \end{pmatrix}  \begin{pmatrix}  0 \\ \eta \end{pmatrix}+o(z) \\
&= -J_h^\perp(\xstar) H_{xx}\Lag(\xstar,\lstar) J_h^{\perp \top}(\xstar) \eta + o(z).
  \end{split} \end{equation}
We notice that the zero dynamics does not depend on $\lambda$.

Finally, by neglecting the high-order terms $o(z)$, the linearization of the zero dynamics is
\begin{equation} \label{linearized_zerodyn}
\dot \eta = -J_h^\perp(\xstar) H_{xx}\Lag(\xstar,\lstar) J_h^{\perp \top}(\xstar) \eta.
\end{equation}
The original nonlinear zero-dynamics is locally asymptotically stable at $\eta=0$ if \eqref{linearized_zerodyn} is asymptotically stable, i.e., if the symmetric matrix $J_h^\perp(\xstar) H_{xx}\Lag(\xstar,\lstar) J_h^{\perp \top}(\xstar)$ is positive definite, which holds if the second-order sufficient conditions reported in Definition \ref{def:2nd} are satisfied. In fact, let $v=J_h^{\perp,\top}(\xstar)w$ for any non-null $w\in\R^{n-m}$; since $J_h(\xstar)J_h^{\perp,\top}(\xstar)=0$ by definition, $J_h(\xstar)v=0$. Then, $v^{\top}H_{xx}\Lag(\xstar,\lstar)v=w^{\top}J_h^{\perp}(\xstar)H_{xx}\Lag(\xstar,\lstar)J_h^{\perp,\top}(\xstar)w>0$.
\end{proof}

\begin{remark}
	In our setting, the linearized zero dynamics in \eqref{linearized_zerodyn} corresponds to the zero dynamics of the linearization of $\plant$, as we prove in the following.
	We notice that the commutativity of the operations of linear approximation and computation of the zero dynamics always holds for single-input, single-output systems; see, e.g., \cite[Remark 4.3.2]{isi95}. However, the commutativity is not guaranteed in general for multiple-input, multiple-output systems. For this reason, it is worth remarking that the class of multiple-input, multiple-output systems considered in this work, characterized by $r=(1,1,\dots,1)\in\R^m$ and $G=-J_H^\top$ according to  \eqref{pezzi}, enjoys the commutativity.

	By Taylor expansion, the linearization of $\plant$ in a neighbourhood of $(\xstar,\lstar)$ is
	\begin{equation}
	 \left\{\begin{array}{l}
	         \dot{x}=-H_{xx}\Lag(\xstar,\lstar)(x-\xstar)-J_h^\top(\xstar)(\lambda-\lstar)\\
	         y= J_h(\xstar)(x-\xstar)
	        \end{array}
\right.
	\end{equation}
Let us consider the mapping $\Phi: \R^n\mapsto \R^n$ as
\begin{equation}
	\begin{pmatrix}
	 \xi\\ \eta
	\end{pmatrix}=
\Phi(x) = \begin{pmatrix}  J_h(\xstar) \\ J_h^\perp(\xstar)  \end{pmatrix}(x-\xstar)
\end{equation}
We notice that $\xi=y$.
In normal form,
\begin{equation}
	\begin{pmatrix}
	 \dot{\xi}\\ \dot{\eta}\end{pmatrix} = \begin{pmatrix}  J_h(\xstar) \\ J_h^\perp(\xstar)   \end{pmatrix}\dot{x}.
\end{equation}
Let us focus on $\eta$, which represents the state variable of the zero dynamics. We have
\begin{equation}
 \begin{split}
  \dot{\eta}&= J_h^\perp(\xstar)\dot{x}\\
  &= -J_h^\perp(\xstar)H_{xx}\Lag(\xstar,\lstar)(x-\xstar).
 \end{split}
\end{equation}
As expected,  $\dot{\eta}$ does not depend on $\lambda$ since $J_h^\perp(\xstar)J_h(\xstar)=0$.
By inversion,
\begin{equation}
\begin{split}
x=& \begin{pmatrix}J_h^{\dagger}(\xstar), & J_h^{\perp,\top}(\xstar)\end{pmatrix}\begin{pmatrix}
    \xi \\ \eta
\end{pmatrix}+\xstar\\
\end{split}
\end{equation}
To study the zero dynamics, we set $\xi=y=0$. Thus,
\begin{equation}
\begin{split}
x-\xstar=& \begin{pmatrix}J_h^{\dagger}(\xstar), & J_h^{\perp,\top}(\xstar)\end{pmatrix}\begin{pmatrix}
    0 \\ \eta
\end{pmatrix}=J_h^{\perp,\top}(\xstar)\eta
\end{split}
\end{equation}
and
\begin{equation}\label{zerodyn}
 \dot{\eta}=-J_h^\perp(\xstar)H_{xx}\Lag(\xstar,\lstar)J_h^{\perp,\top}(\xstar)\eta
\end{equation}
which is equal to \eqref{linearized_zerodyn}.
\end{remark}

\subsection{Global exponential convergence of FL-CMO for strongly convex problems}\label{sub:FL_conv_stronglyconvex}

In this section, we study the global convergence of FL-CMO dynamics, defined by \eqref{eq:flcmo}, for strongly convex problems. More precisely, we consider the same setting in which we analyse the convergence of PI-CMO and described by assumptions \ref{ass:sc} and \ref{ass:aff}, namely $f(x)$ in \eqref{opt} is strongly convex and twice differentiable; $h(x)=Cx+d$ where $C\in\R^{m,n}$ is full rank and $\alpha_1 I_m \preceq CC^\top \preceq \alpha_2 I_m$.
Under these assumptions, $J_h(x) = C$, and the control law for \eqref{uncontrolled} defined by \eqref{decoupling_control}, \eqref{Abr1} is:
\begin{equation}
\lambda(t) = -(CC^\top)^{-1}(C\nabla f(x(t))+v(t))
\end{equation}
where $v(t)=\mathcal{G}(y(t))$, see \eqref{linear_control}, and $y(t)=Cx(t)+d$.

Moreover, according to \cite[Lemma 1]{qu19}, there exists a symmetric $B(x)\in\R^{n,n}$ satisfying $\beta_1 I \preceq B(x) \preceq \beta_2 I$ for some $0< \beta_1<\beta_2$ such that \eqref{lem:qu} holds. As in Sec. \ref{sec:PI}, for brevity we use the notation $B=B(x)$.

\begin{theorem}[Global exponential convergence of FL-CMO]\label{th:glob_exp_fl}
Let assumptions \ref{ass:sc} and \ref{ass:aff} hold.
Let us choose $\mathcal{G}$ in \eqref{linear_control} such that $\bar y=0$ is a globally exponentially stable equilibrium for $\dot y = \mathcal{G}(y)$, with convergence rate $\mu_g\geq \beta_1$, i.e., there exists a constant $c_g \in \R_+$ such that
    \begin{equation}\label{goodG}
        \norm{y(t)}_{2} \leq c_g e^{-\mu_g t}.
    \end{equation}
    Then, the global optimum $\xstar$ of problem \eqref{opt} is a globally exponentially stable equilibrium for FL-CMO with rate $\beta_1$, i.e., there exists a constant $c\in \R_+$ such that
    \begin{equation}\label{FLCMO_thesis}
        \lVert x-\xstar \rVert_2 \leq c e^{-\beta_1 t}.
    \end{equation}
\end{theorem}

\begin{proof}
Under assumptions \ref{ass:sc}-\ref{ass:aff}, the global change of coordinates defined by \eqref{def:z} and \eqref{def:phi} is
    \begin{equation}
        z = \Phi(x)= \begin{pmatrix} \zeta \\ \eta \end{pmatrix} = \begin{pmatrix}  Cx+d \\ C^\perp (x-\xstar) \end{pmatrix}
        \label{eq:linear_phi}
    \end{equation}
    where the rows of $C^\perp\in\R^{m-n,n}$ are an orthonormal basis for the null space of the rows of $C$.

    Since \eqref{eq:linear_phi} is an affine transformation,
    \begin{equation}\label{invert}
        x-\xstar = \begin{pmatrix}  C \\ C^\perp \end{pmatrix}^{-1} (z-\zstar)
    \end{equation}
    where $\zstar=\Phi(\xstar)$.

    Given the spectral norm  $\sigma:=\norm{\begin{pmatrix}  C \\ C^\perp \end{pmatrix}^{-1}}_{2}$ and by applying the triangle inequality we obtain
    \begin{equation}\label{bound_x}
         \norm{x-x^\star}_2 \leq \sigma \norm{z-z^\star}_2
         \leq \sigma \left(\norm{\begin{pmatrix}
            \zeta-\zeta^\star \end{pmatrix}}_2+
\norm{\begin{pmatrix}
             \eta-\eta^\star \end{pmatrix}}_2\right).
   \end{equation}
Then, we study the convergence of $\norm{x-x^\star}_2$ based on the convergence of $\norm{\begin{pmatrix}
            \zeta-\zeta^\star\end{pmatrix}}_2$ and $\norm{\begin{pmatrix}
            \eta-\eta^\star \end{pmatrix}}_2$.

Since $\zeta=y$ and $\zeta^\star=0$, the global exponential convergence of $\zeta$ is given by \eqref{goodG}.

As to $\eta$, by using $\nabla_x \mathcal{L}(\xstar,\lambda(\xstar))=0$, the zero dynamics is
\begin{equation}\label{eq:lin_zerodyn1}\begin{aligned}
        \dot \eta &= C^\perp \dot x = C^\perp \left[\nabla_x \mathcal{L}(x,\lambda(x))-\nabla_x \mathcal{L}(x^\star,\lambda(x^\star))\right] = \\
        &=-C^\perp \left[I+C^\top(CC^\top)^{-1}C) \big(\nabla f(x)-\nabla f(x^\star)\big)\right]\\
        &=-C^\perp \big(\nabla f(x)-\nabla f(x^\star)\big)
    \end{aligned}\end{equation}
    where the last equality follows from $\mathcal{G}(0)=0$ and $C^\perp C^\top=0$.

According to \cite[Lemma 1]{qu19} on strongly convex functions, see \eqref{lem:qu}, we can write
$\dot \eta = C^\perp B [x-\xstar]$.
Moreover,
    \begin{equation*}\begin{aligned}
        x-x^\star &= C^\top(CC^\top)^{-1} (\zeta-\zeta^\star) + C^{\perp \top} (\eta-\eta^\star) = C^{\perp \top} \eta
    \end{aligned}\label{eq:diffxaseta_flcmocvx}\end{equation*}
    where we use  $\zeta=\zeta^\star=0$ by definition of zero dynamics  and  $\eta^\star=0$.
    In conclusion,
    \begin{equation}
        \dot \eta = -C^\perp B C^{\perp \top} \eta.
    \end{equation}

    From \cite[Lemma 1]{qu19}, $B- \beta_1 I\succeq 0$. Moreover, by applying \cite[Observation 7.1.8]{hor13}, we have
    \begin{equation}
    C^\perp (B-\beta_1 I) C^{\perp \top}=C^\perp BC^{\perp \top} -\beta_1 I \succeq 0.
    \end{equation}
    
    Thus, the matrix $C^\perp BC^{\perp \top}$ is Hurwitz and
    \begin{equation}\label{rate_eta}
        \norm{\eta} \leq c_\eta e^{-\beta_1 t}.
    \end{equation}
   By merging  \eqref{bound_x}, \eqref{goodG}, and \eqref{rate_eta}, we get
    \begin{equation}
        \norm{x-\xstar}_2 \leq \sigma \left(  c_g e^{-\mu_g t} + c_\eta e^{-\beta_1 t} \right) \leq  \sigma (  c_g  + c_\eta) e^{-\beta_1 t}
    \end{equation}
    which proves \eqref{FLCMO_thesis} with $c=\sigma (c_g+c_\eta)$.
\end{proof}

\begin{remark}\label{rem:faster} Theorem \ref{th:glob_exp_fl} shows that FL-CMO enjoys a better convergence rate than PI-CMO in strongly convex problems. In fact, the rate $\beta_1$ of FL-CMO is always larger than $\frac{1}{2}\mu< \beta_1$ evaluated in Theorem \ref{theo:pi} for PI-CMO.
\end{remark}

\begin{remark}The convergence rate of FL-CMO uniquely depends on $\beta_1$, i.e., on $f(x)$, while it is independent from the constraints $h(x)$. This is not the case for PI-CMO and PDGD which depend on constants $\alpha_1$ and $\alpha_2$, that are related to the constraints.
\end{remark}

\section{Numerical results}\label{sec:NR}
In this section, we illustrate five numerical examples that validate the theoretical convergence results and prove the effectiveness of the proposed approaches PI-CMO and FL-CMO compared to state-of-the-art algorithms.

In Sec. \ref{sub:first}, we test PI-CMO and FL-CMO  in a convex problem and analyze their convergence speed.
Then, in Sec. \ref{sub:second}, we investigate the effectiveness of PI-CMO in a non-convex problem in which FL-CMO is not feasible.

In the succeeding two examples, we test FL-CMO  in non-convex problems, namely a gray-box system identification problem in Sec. \ref{sub:third} and a large-scale real-world chemical problem in Sec. \ref{sub:fourth}.

Finally, in Sec.\ref{sub:fifth}, we illustrate the behavior of different instances of FL-CMO  based on different choices of $\mathcal{G}$.

\subsection{PI-CMO and FL-CMO vs PDGD in convex optimization}\label{sub:first}
In the first example, we resort to the quadratic optimization problem with linear constraints proposed in \cite[Section IV.A]{qu19}.
Specifically, we consider
\begin{equation}
\begin{split}
&\min_{x\in\R^n} \frac{1}{2}x^\top W x\\
&\text{s.t.} \quad Cx+d=0
\end{split}
\end{equation}
where $W=10I+W_0W_0^\top \in \R^{n,n}$ is positive definite; $W_0\in\R^{n,n}$, $C\in\R^{m,n}$ and $d\in\R^m$ have independent and normally distributed components.

To solve this problem, we implement and compare PI-CMO, FL-CMO and PDGD \cite{qu19}.
Since the cost function is quadratic and the constraints are linear, the  PI-CMO dynamics is linear time-invariant, as illustrated in Sec. \ref{sub:quad}. 

Corollary \ref{cor:PI} shows that the convergence rate bound of PI-CMO compares favorably with the one of PDGD for strongly convex problems with linear constraints. Moreover, as noticed in Remark \ref{rem:faster}, the convergence rate bound of FL-CMO enhances that of PI-CMO.
In this section, we support these theoretical results with a numerical analysis of the convergence speed.

We consider $n=50$ variables and $m\in [2,26]$ constraints. For PDGD, we set $\eta=\frac{\beta_1 \beta_2}{\alpha_2}$ to obtain the best convergence rate bound; see \eqref{mupdgd}. For PI-CMO, we set $K_i$ and $K_p$ to meet the conditions of Corollary \ref{cor:PI}, specifically  $K_i=\frac{50}{\epsilon}\frac{\beta_1\beta_2}{8\alpha_2}$, $K_p=\epsilon\frac{2K_i}{\beta_2}$ where $\epsilon=\frac{4}{5}\left(1-\frac{\alpha_1}{8\alpha_2}\right)$.
For FL-CMO, we design a static $\mathcal{G}$ with convergence rate $\mu_g=5\beta_1$.

For all the algorithms, we run the Euler discretized versions with the discretization step size that guarantees stability, see \cite[Section III.C]{qu19} for details.
We perform 400 random runs.
\begin{figure}
\centering
 \includegraphics[width=0.9\linewidth]{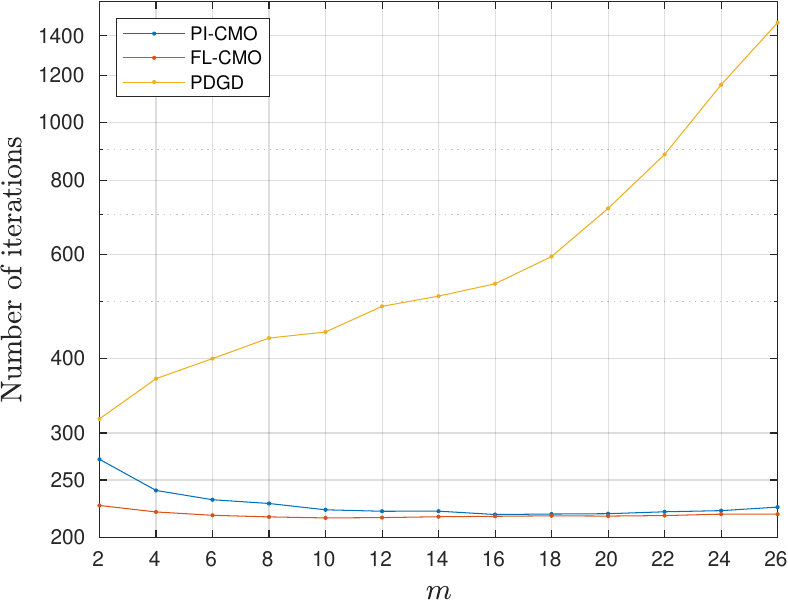}
 \caption{Number of iterations to converge for PI-CMO, FL-CMO, and PDGD in a convex quadratic example, with $m$ constraints. The results are averaged over 400 random runs.}
 \label{fig:PI_vs_PDGD}
\end{figure}
Fig.~\ref{fig:PI_vs_PDGD} shows the average number of iterations required to converge to the desired minimum against $m$. For all the considered values of $m$, PI-CMO requires fewer iterations on average than PDGD, and FL-CMO is even faster than PI-CMO, as expected from the theoretical results on the convergence rate presented in Sec.~\ref{sub:PI_conv} and \ref{sub:FL_conv_stronglyconvex}.

\begin{remark}\label{rem:PI_vs_CMO}
FL-CMO enjoys a better convergence rate than PI-CMO. Nevertheless, FL-CMO presents two main computational drawbacks that make PI-CMO more effective in some contexts.
The first drawback is that FL-CMO requires the inversion of the $m\times m$ matrix $J_h(x)J_h(x)^\top$ in \eqref{decoupling_control} at each iterative step, which generates a computational bottleneck.
If the constraints are affine,  $J_h(x)J_h(x)^\top=CC^\top$ is constant, which basically solves the issue: the inversion is computed once. In the other cases, we can exploit the structure of $A(x)$ to perform a suitable factorization; see Remark~\ref{rem:bottleneck}.
Another possible workaround is to consider the approximate computation of the inverse; see Sec. \ref{sub:fifth}.
When the inversion is not feasible, PI-CMO is a good alternative to FL-CMO.
The second drawback is that implementing FL-CMO requires the invertibility of $J_h(x)J_h(x)^\top$. Therefore, this method cannot be implemented if $m>n$, which is not required for implementing PI-CMO, as we show in the next numerical example.
\end{remark}

\subsection{Shidoku puzzle}\label{sub:second}
\begin{figure}[h]
    \centering
    \includegraphics[width=0.35\linewidth]{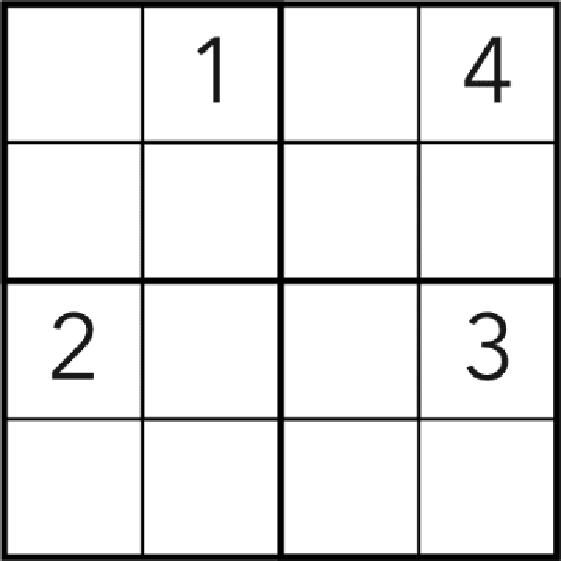}
    \caption{Shidoku puzzle to be solved.}
    \label{fig:shidoku}
\end{figure}

Shidoku is a 4x4 version of the popular 9x9 Sudoku puzzle. Given an initial scheme, as reported in Fig. \ref{fig:shidoku}, the aim is to fill the empty cells with integers $x_{i,j} \in \{1,2,3,4\}$ such that each row, each column, and each 2x2 corner block contains the integers $1,2,3,4$.
We can formulate the game in terms of the solution of polynomial equations on the values $x_{i,j}$ of each cell $(i,j)$, $i,j=1,\dots, 4$. First of all, $x_{i,j}\in\{1,2,3,4\}$ is guaranteed by $\prod_{h=1}^4 (x_{i,j} - h) = 0$. Then, we obtain no repetition in groups of 4 cells by imposing the product equal to $24$ and the sum equal to $10$.
In summary, given the corner blocks
\begin{equation*}
    \begin{split}
    B_1 &= \{(1,1), (1,2), (2,1), (2,2)\}\\
    B_2 &= \{(1,3), (1,4), (2,3), (2,4)\}\\
    B_3 &= \{(3,1), (3,2), (4,1), (4,2)\}\\
    B_4 &= \{(3,3), (3,4), (4,3), (4,4)\}.
\end{split}
\end{equation*}
we have
\begin{equation}\label{eq:shidoku} \begin{aligned}
    &\text{Columns: for } j = 1,\dots,4,\\
    &~~~~\sum_{i=0}^4 x_{ij} = 10, \quad \prod_{i=0}^4 x_{ij} = 24. \\
    &\text{Rows: for } i = 1,\dots,4,\\
    &~~~~\sum_{j=0}^4 x_{ij} = 10, \quad \prod_{j=0}^4 x_{ij} = 24.\\
    &\text{Blocks: for }  k = 1,\dots,4,\\
    &~~~~\sum_{(i,j)\in B_k} x_{ij} = 10, \quad \prod_{(i,j)\in B_k} x_{ij} = 24.\\
    &\text{Restriction to integers: for } i,j = 1,\dots,4,\\
    &~~~~\prod_{h=1}^4 (x_{ij} - h) = 0. \\
    &\text{Initial conditions as in Fig. \ref{fig:shidoku}:}\\
    &~~~~x_{1,2} = 1 \quad x_{1,4} = 4 \quad x_{3,1} = 2 \quad x_{3,4} = 3.
\end{aligned}\end{equation}
Equations \eqref{eq:shidoku} represent non-convex polynomial constraints. We can solve the corresponding optimization problem through the proposed CMO framework if we associate them with any constant cost function. In particular, we use PI-CMO to test its convergence and effectiveness in non-convex problems. We cannot apply FL-CMO since we have $m=40$ equations in $n=12$ variables, and $m>n$ is not consistent with Assumption \ref{ass:fl2}.

As to PI-CMO, we set $K_i = 1, K_p = 0.1$ and we generate random initial conditions according to $x_{i,j}(0) = \vert \xi_{i,j} \vert, \xi_{i,j} \sim \mathcal{N}(0,1)$, for each $i,j=1,\dots, 4$ and $\lambda_k(0) \sim \mathcal{N}(0,1), k = 1,\dots,m$, where $\mathcal{N}(0,\sigma)$ denotes the normal distribution with zero mean and variance $\sigma^2$. We integrate the ordinary differential equations that describe the closed-loop dynamics thanks to the MATLAB \rm{ode45} solver in the time interval $[0,100]$ seconds, corresponding to approximately $151700$ iterations.
We perform 20 runs with different random initial conditions. PI-CMO is always convergent in the given time interval.

\begin{figure}
    \centering
    \includegraphics[width=0.95\linewidth]{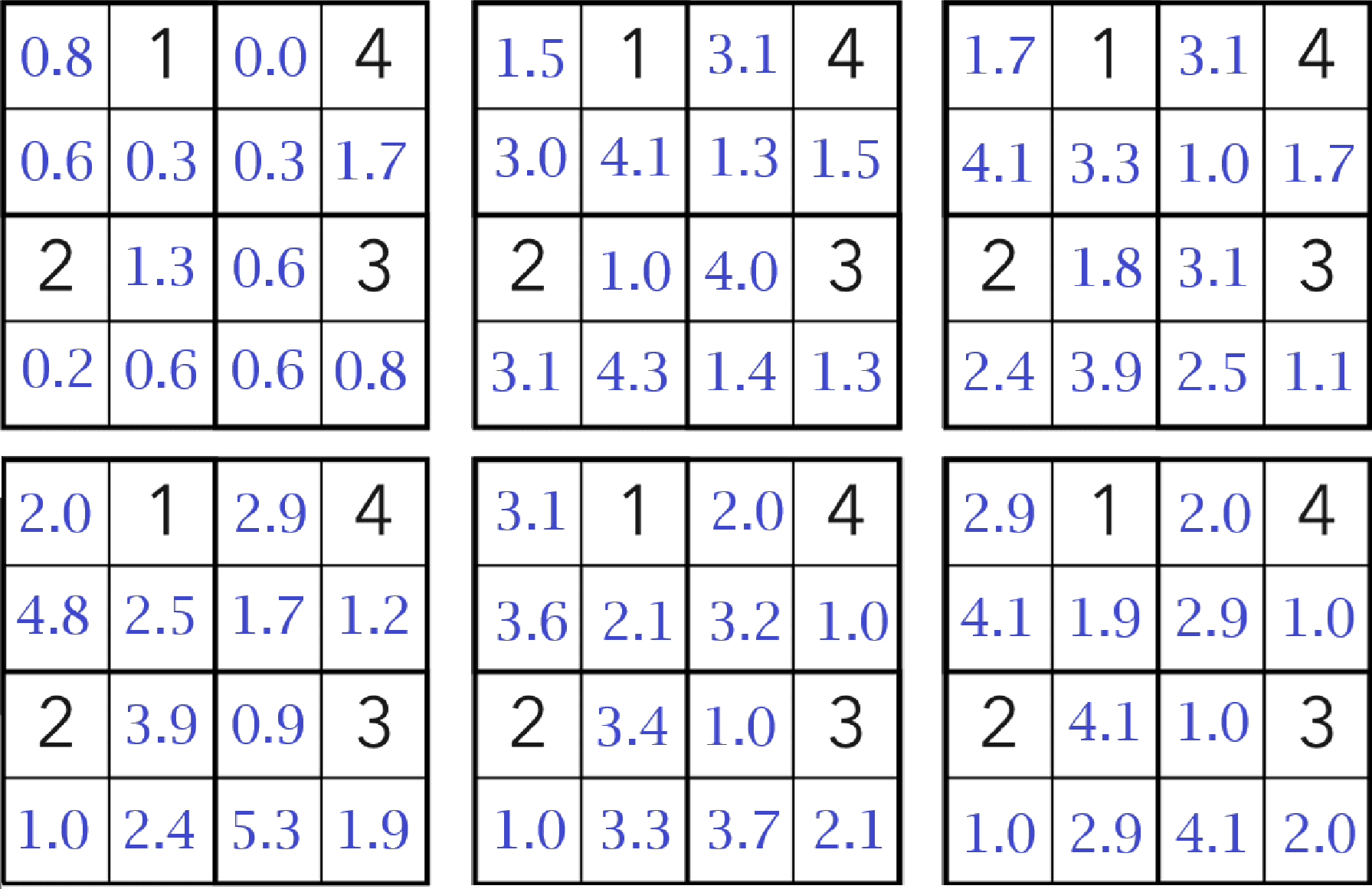}
    \caption{PI-CMO solves Shidoku puzzle. Evolution of the solution of the optimization variables at six equispaced sampling instants between initialization and convergence step $91027$, from top left to bottom right.}
    \label{fig:shidoku_sol_evol}
\end{figure}

\begin{figure}
    \centering
    \includegraphics[width=0.35\linewidth]{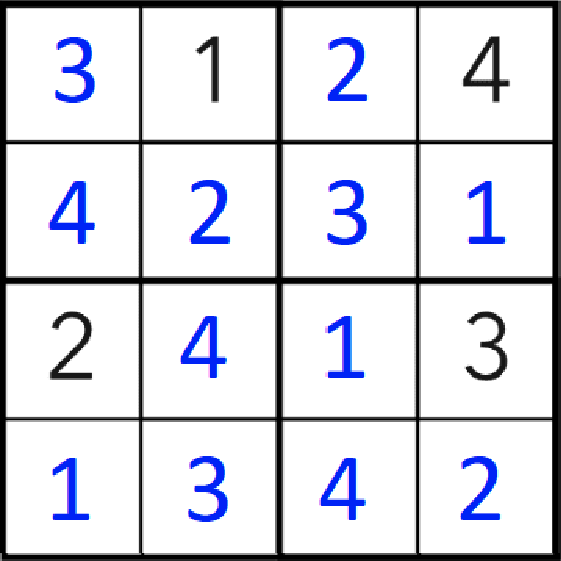}
    \caption{PI-CMO solves Shidoku puzzle: final correct solution.}
    \label{fig:shidoku_sol}
\end{figure}
We show an instance in Fig. \ref{fig:shidoku_sol_evol}, where we depict the evolution of the optimization variables for six equispaced sampling instants between the random initialization and the convergence to the correct solution at iteration 91027, shown in Fig. \ref{fig:shidoku_sol}.

In conclusion, this test shows that PI-CMO is convergent even in a non-convex problem that does not satisfy assumptions \ref{ass:sc} and \ref{ass:aff}.
For further investigation, we compare PI-CMO to two state-of-the-art approaches for non-convex constrained optimization, interior-point method (IPM) and sequential quadratic programming, through the \rm{fmincon} function in MATLAB. We perform 20 runs with random initial conditions for the two of them.
We observe that IPM fails in all the runs due to numerical issues. More precisely, the linear system of KKT conditions to solve at each iteration is poorly conditioned, which affects the solution; see, e.g., \cite[Chapter 19]{nocedal2006numerical} for details. On the other hand, sequential quadratic programming converges to an infeasible point in all the runs.

\subsection{Gray-box non-linear system identification}\label{sub:third}
In this example, we consider the identification of the discrete-time nonlinear system described by the following regressor form
\begin{equation}\begin{split}
    y(k) = & \theta_1 e^{-(y(k-1))^2} +\theta_2 (u(k-1))^2 + \\ \qquad & \theta_3 u(k-2)y(k-1) + \theta_4 (u(k-2))^{\theta_5}
\end{split}\end{equation}
using $N = 400$  measurements of the input sequence $u(k)$ and noisy measurements of the output sequence $\tilde{y}(k) = y(k) + \zeta(k)$, where $\zeta(k)$ represents an unknown measurement noise.

Since the system is not affine in the parameters $\theta$, we cannot employ standard solutions based on least-squares regression.
We look for the parameter vector $\theta^\star\in\R^5$ that minimizes the $\ell_2$ norm of $\zeta = [\zeta(1),\dots,\zeta(N)]^\top$, i.e., the energy of the noise sequence $\zeta(k)$.
Introducing additional optimization variables $y \in \R^N$ representing the noise-free output samples, the considered identification problem is equivalently rewritten as the following non-convex constrained problem:
\begin{equation} \begin{split}
        [\theta^\star,y^\star] & = \argmin{\theta \in \R^5, y \in \R^N } \sum_{k=1}^N (y_k - \tilde{y}(k))^2 \\
                       & \text{subject to:}                                                                  \\
                       & -y_k + \theta_1 e^{-y_{k-1}^2} +\theta_2 (u(k-1))^2 +                             \\ & \theta_3 u(k-2)y_{k-1} + \theta_4 (u(k-2))^{\theta_5}=0 \\ & \qquad k = 3,\dots,N.
        \label{probl_fl_example}
    \end{split} \end{equation}

We solve problem \eqref{probl_fl_example} using the feedback linearization method developed in Section \ref{sec:FL}. We define the controller $\mathcal{G}$ in  \eqref{linear_control} according to \eqref{linear_control_details} , with $K_i=1$ for each $i=1,\dots,m$. 

We integrate the FL-CMO ordinary differential equations in the time interval $t = [0,20]$ seconds using Euler discretization with step size $10^{-2}$ seconds. We set normally distributed initial conditions.

\begin{figure}
    \centering
    \input{esNNARX_state_evol}
    \caption{Gray-box non-linear system identification: evolution of the estimation error. We denote by $\hat \theta(k)$ the estimates obtained with the feedback linearization-based algorithm at iteration $k$, and by $\theta_{\rm true}$ the true parameter vector.}
     \label{fig:sysid_param}
\end{figure}

\begin{figure}
    \centering
    \input{esNNARX_cns_evol}
    \caption{Gray-box non-linear system identification: evolution of the $\ell_\infty$-norm of the constraints function $h$. In the figure, $x(k)$ refers to the estimate, at iteration $k$, of $y$ and $\theta$ in \eqref{probl_fl_example}.}
    \label{fig:sysid_cns}
\end{figure}

Fig. \ref{fig:sysid_param} shows the $\ell_2$-norm of the estimation error as a function of the algorithm iteration. Fig. \ref{fig:sysid_cns} shows the evolution of the $\ell_\infty$-norm of the constraints: as expected, it converges to zero.

For comparison, we solve problem \eqref{probl_fl_example} through IPM implemented with the \rm{fmincon} MATLAB function and initialized with the same initial conditions used for the feedback linearization method. We use MATLAB R2021b on a processor i7-10700, 2.90GHz with 32 GB of DDR4 RAM.
\begin{table}
\caption{Gray-box non-linear system identification: comparison of parameter estimates for the proposed feedback linearization approach and IPM, denoted by $\hat \theta$ and $\theta_{\rm IPM}$, respectively.}
    \label{tab:comparison}
    \centering
    \begin{tabular}{c | c  c }
        $\theta_{\rm true}$  &  $\hat \theta$ & $\theta_{\rm IPM}$ \\
       \hline
            $0.5$    &      $0.496$    &    $ 0.489 $  \\
           $-0.3$   &     $ -0.300 $    &    $  -0.300$   \\
           $ -0.7 $   &   $  -0.699$   &     $ -0.695$   \\
           $-0.35 $    &   $ -0.345$   &     $ -0.338$  \\
             $0.8  $  &   $   0.815 $   &      $ 0.804 $
   \end{tabular}
\end{table}
Table \ref{tab:comparison} compares the estimated parameters. We observe that both feedback linearization and IPM provide accurate estimates, the respective errors being $\|\theta_{\rm true}-\hat \theta\|_2=0.0163$ and $\|\theta_{\rm true}-\theta_{\rm IPM}\|_2=0.0175$.

Let us analyze the computational complexity. The time required by FL-CMO is $9.8$ seconds, while IPM  requires $64.4$ seconds. Moreover, the proposed feedback linearization method is less memory expensive compared to IPM, as expected when we compare first-order and second-order methods because it requires storing only the Jacobian of the constraints and not the Hessian matrix.
More precisely, given $m\leq n$, the FL-CMO requires $n+nm+\frac{1}{2}m^2+m$ floating point numbers, while IPM requires $n+n^2+nm+3m$ floating point numbers. In both cases, we assess this number by assuming to store only the triangular part of the symmetric matrices.
Concerning IPM, the leading term $n^2$ is due to the sum of $m+1$ Hessian matrices, one for $f(x)$ and $m$ for $h_i(x)$, of dimension $n \times n$. Thus, we need $\frac{1}{2}n^2$ variables to compute the current Hessian, and a further $\frac{1}{2}n^2$ ones are used to store the partial sum.

In conclusion, there is a gain of $2m+n^2-\frac{1}{2}m^2$ floating-point numbers using FL-CMO instead of IPM.

Moreover, each iteration in feedback linearization requires $O(m^3)$ floating-point operations (FLOPs) to invert the matrix $J_h J_h^\top \in \real^{m\times m}$. In contrast, each iteration of IPM requires $O((n+m)^3)$ FLOPs to solve a linear system of dimension $n+m$.

We notice that similar considerations apply to sequential quadratic programming, which needs approximately the same memory and computations as IPM; see, e.g., \cite[Chapter 18]{nocedal2006numerical}.

\subsection{Industrial chemical process problem}\label{sub:fourth}
We consider a problem that arises in the context of industrial chemical processes: the propane, isobutane, and n-butane nonsharp separation presented in \cite{agg90}. It is about a three-component feed mixture required to separate products into two three-component products. This problem is included in the benchmark suite of real-world, non-convex problems proposed in \cite{kum20} to test optimization algorithms.
The mathematical formulation of the problem, reported in \cite[Sec. 2.1.5]{kum20}, is as follows. Given $x=(x_1,\dots,x_{48})$, minimize $f(x)$ defined as
\begin{equation*}
\begin{split}
f(x) &= c_{11}+\left(c_{21}+c_{31} x_{24}+c_{41} x_{28}+c_{51} x_{33}+c_{61} x_{34}\right) x_{5}\\&+c_{12}+\left(c_{22}+c_{32} x_{26}+c_{42} x_{31}+c_{52} x_{38}+c_{62} x_{39}\right) x_{13}
\end{split}
\end{equation*}
where
$c_{11}=0.23947$, $c_{12}=0.75835$, $c_{21}=-0.0139904$, $c_{22}=-0.0661588$, $c_{31}=0.0093514$, $c_{32}=0.0338147$, $c_{41}=0.0077308$, $c_{42}=0.0373349$, $c_{51}=-0.0005719$, $c_{52}=0.0016371$, $c_{61}= 0.0042656 $, $c_{62}=0.0288996$,
subject to
\begin{equation*}
\begin{aligned}
    &x_{4}+x_{3}+x_{2}+x_{1}=300,\quad &&  x_{6}-x_{8}-x_{7}=0,\\
    &x_{9}-x_{12}-x_{10}-x_{11}=0,\quad &&x_{14}-x_{17}-x_{15}-x_{16}=0,\\
    &x_{18}-x_{20}-x_{19}=0,\quad &&x_{6} x_{21}-x_{24} x_{25}=0, \\
    &x_{14} x_{22}-x_{26} x_{27}=0,\quad && x_{9} x_{23}-x_{28} x_{29}=0, \\
    &x_{18} x_{30}-x_{31} x_{32}=0,\quad && x_{25}-x_{5} x_{33}=0, \\
    &x_{35}-x_{5} x_{36}=0,\quad &&x_{37}-x_{13} x_{38}=0, \\
    &x_{27}-x_{13} x_{39}=0,\quad &&x_{32}-x_{13} x_{40}=0,\\
    &x_{25}-x_{6} x_{21}-x_{9} x_{41}=0,\quad&& x_{29}-x_{6} x_{42}-x_{9} x_{23}=0, \\
    &x_{35}-x_{6} x_{43}-x_{9} x_{44}=0,\quad &&x_{37}-x_{14} x_{45}-x_{18} x_{46}=0,\\
    &x_{27}-x_{14} x_{22}-x_{18} x_{47}=0,\quad &&x_{32}-x_{14}x_{48}-x_{18}x_{30}=0,\\
    &0.33 x_{1}+x_{15} x_{45}-x_{25}=0,\quad &&0.33 x_{1}+x_{15} x_{22}-x_{29}=0, \\
    &0.33 x_{1}+x_{15} x_{48}-x_{35}=0,\quad && 0.33 x_{2}+x_{10} x_{41}-x_{37}=0, \\
    &0.33 x_{2}+x_{10} x_{23}-x_{27}=0,\quad &&0.33 x_{2}+x_{10} x_{44}-x_{32}=0, \\
    & x_{33}+x_{34}+x_{36}=1,\quad && x_{21}+x_{42}+x_{43}=1,\\
    & x_{41}+x_{23}+x_{44}=1,\quad && x_{38}+x_{39}+x_{40}=1, \\
    & x_{45}+x_{22}+x_{48}=1,\quad && x_{46}+x_{47}+x_{30}=1, \\
    & x_{43}=0, \quad && x_{46}=0,
\end{aligned}
\end{equation*}
\begin{equation*}
\begin{aligned}
    &0.33 x_{3}+x_{7} x_{21}+x_{11} x_{41}+x_{16} x_{45}+x_{19} x_{46}=30,\\
    &0.33 x_{3}+x_{7} x_{42}+x_{11} x_{23}+x_{16} x_{22}+x_{19} x_{47}=50,\\
    &0.33 x_{3}+x_{7} x_{43}+x_{11} x_{44}+x_{16} x_{48}+x_{19} x_{30}=30, \\
\end{aligned}
\end{equation*}
with bounds
\[
\begin{aligned}
& 0 \leq x_{1}, \ldots, x_{20} \leq 150,\\
& 0 \leq x_{25}, x_{27}, x_{32}, x_{35}, x_{37}, x_{29} \leq 30, \\
& 0 \leq x_{21}, x_{22}, x_{23}, x_{30}, x_{33}, x_{34}, x_{36}, x_{38}, x_{39} \leq 1, \\
& 0 \leq  x_{40}, x_{42}, x_{43}, x_{44}, x_{45},x_{46}, x_{47}, x_{48} \leq 1,\\
&0.85 \leq x_{24}, x_{26}, x_{28}, x_{31} \leq 1.
\end{aligned}
\]    % two constraints per line

The problem consists of $n=48$ optimization variables, a bilinear cost function, and $m=38$ linear and bilinear equality constraints. The overall problem is non-convex. Moreover, there are $47$  constrained variables in bounded intervals, which give rise to $94$ inequality constraints. To deal with the inequality constraints, we reformulate them by using squared-slack variables, as discussed, e.g., in \cite[Sec. 3.3.2]{ber99}. The basic idea is that we can rewrite any inequality  $g(x)\leq 0$ as $g(x)+z^2=0$, where $z\in\R$ is a slack variable.

For this problem, the benchmark objective value is  $f(x)=2.1158$, as reported in  \cite[Table 3]{kum20}.

We perform the optimization through FL-CMO, implemented in MATLAB R2021b, on a processor i7-10700, 2.90GHz with 32 GB of DDR4 RAM.
As in the previous example, we design the $\mathcal{G}$ controller in Eq.~\eqref{linear_control} by pole placement. We place the closed-loop pole at $-10$. We integrate the closed-loop differential equations in the time interval $[0,3]$ seconds using Euler discretization with step size $5 \cdot 10^{-5}$ seconds.

\begin{figure}
    \centering
    \includegraphics[width=0.95\linewidth]{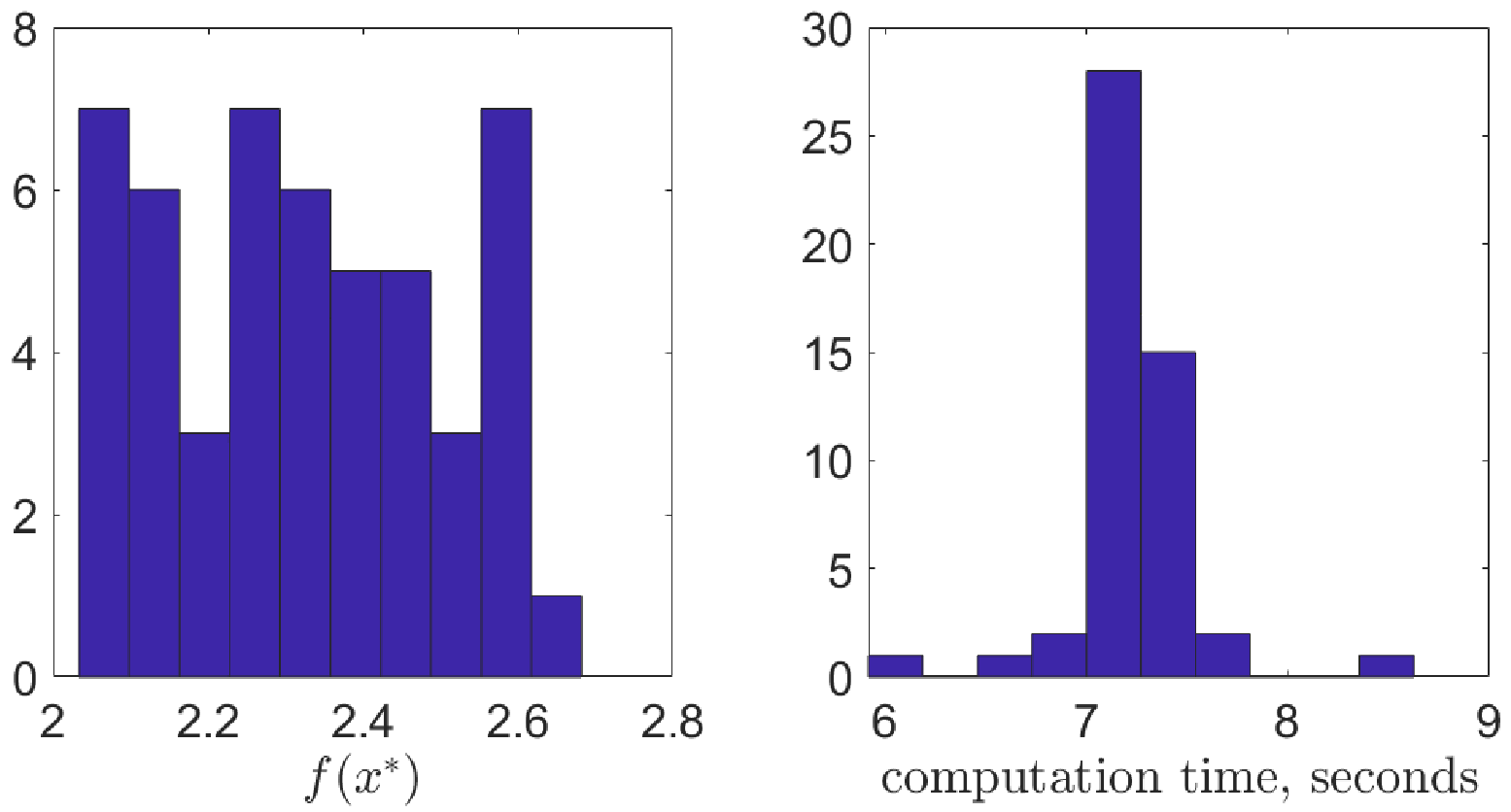}
    \caption{Industrial chemical process problem: distribution of the achieved objective values (left) and distribution of the required time to converge (right)}
    \label{fig:chem}
\end{figure}

We randomly set the initial conditions with uniform distribution in $[0,50]$. We run the algorithm $50$ times with different realizations of the initial conditions. In all the runs, the algorithm converges to feasible solutions within a tolerance of $10^{-7}$ for each constraint.  The achieved objective values are distributed according to the histogram in Fig. \ref{fig:chem}; the mean value is $2.332$, and the standard deviation is $0.1852$. Moreover, the best achieved objective value is $2.0333$, which improves the benchmark $2.1158$. In Fig. \ref{fig:chem}, we also report the distribution of the time required to converge; the mean value is $7.22$ seconds, and the standard deviation is $0.33$ seconds.

We finally remark that IPM is an alternative solver for this problem, but it suffers the increase of required memory to store the Hessian matrix; see the computational complexity analysis in Section \ref{sub:third}.

In conclusion, this experiment proves that the proposed feedback linearization approach is valuable in large-scale, non-convex optimization problems arising from real-world applications. Moreover, using squared-slack variables is a feasible approach to deal with inequality constraints in this example, even if it increases the number of optimization variables and introduces additional non-convex constraints.

\subsection{Robust optimization in the presence of inexact data}\label{sub:fifth}
\label{sec:example_robustness}
As noticed in Sec. \ref{sec:FL}, FL-CMO requires the inversion of  $A(x)\in \R^{m\times m}$, which is a crucial point from different perspectives. On the one hand, the inverse of a poorly conditioned matrix is sensitive to small perturbations in the data. On the other hand, inversion is computationally intense, and one can use approximate solutions to reduce the burden. We deal with an inexact matrix $A(x)^{-1}$ in both cases. 

In the following example, we demonstrate how the flexibility in the design of the controller $\mathcal{G}$ for FL-CMO leads to algorithms characterized by different robustness properties in the presence of such inaccurate data.
Moreover, we compare PI-CMO and FL-CMO equipped with a PI control law for $v(t)$.

We consider the quadratic problem
\begin{subequations}
    \begin{align}
    \min_{x \in \R^3}\, &\frac{1}{2}x^\top \begin{pmatrix}
        1& 0& 1\\ 0& 4 &-2\\ 1& -2& 8
    \end{pmatrix} x + \begin{pmatrix}
        -1\\ 2\\ -1
    \end{pmatrix}^\top x \\
    & \textrm{s.t.} \quad x_1 + 1=0,\\
    & \qquad \,  3x_1+ 2x_2 -4x_3=0. 
    \end{align}
\end{subequations}

Then $A(x)^{-1}=-\left(C C^\top\right)^{-1}=\begin{pmatrix}-1.45&0.15\\0.15&-0.05\end{pmatrix}$. Based on previous considerations, we assume to know only a perturbed version $A(x)^{-1}\odot (\textbf{1}_2+P)$ where $\textbf{1}_2$ is the $2 \times 2$ matrix of all ones, $\odot$ is the Hadamard product, and $P \in \R^{2,2}$ is such that each of its entries $P_{ij}$ is bounded by $\vert P_{ij} \vert \leq 0.25$, for $1 \leq i,j \leq 2$. % 

For FL-CMO, we consider the static controller $v(t)=\mathcal{G}(y(t)) = Ky(t)$ with $K=-4$, which retraces the algorithm proposed in \cite{all24}, and the dynamic PI controller
\begin{equation}\label{eq:exampleperturb_piG}
    v(t)=\mathcal{G}(y(t)) = k_p y(t) + k_i \int_{0}^t y(\tau) \textrm{d}\tau.
\end{equation}
where $k_i= -1$ and $k_p = -4$.
Additionally, we implement FL-CMO with exact $A(x)^{-1}$ as a benchmark. Finally, we also consider PI-CMO with the same coefficients $k_i$ and $k_p$ for further comparison. 

We perform $40$ simulations of the four dynamics, where for each run, we start from different initial conditions and sample different $P_{ij}$ from the uniform distribution in $[-0.25,0.25]$. We show the simulation results in terms of the average distance from the optimum $\norm{x(t)-\xstar}_2$ and the maximum constraint violation $\norm{h(x(t))}_\infty$ in Fig. \ref{fig:xt_perturbed_example} and Fig. \ref{fig:xt_perturbed_example2}, respectively. As expected, the exact FL-CMO and PI-CMO converge to the correct solution.
Interestingly, we can see that the pole placement controller makes FL-CMO converge quickly to a wrong solution under the effect of perturbation. On the other hand, by introducing the integral term in $\mathcal{G}$, the inexact FL-CMO converges to the correct solution more quickly than PI-CMO.

In conclusion, this example suggests that FL-CMO with a dynamic controller is more robust to data perturbation and inaccuracies. Future work will analyze the robustness properties of FL-CMO in more general settings.

\begin{figure}
    \input{robexample_xt_avg_norm}
    \caption{Perturbed quadratic problem: optimization variables' trajectories $x(t)$ for the different algorithms. Exact FL-CMO refers to the unperturbed case.}
    \label{fig:xt_perturbed_example}
\end{figure}

 \begin{figure}
    \input{robexample_ht_avg_norm}
    \caption{Perturbed quadratic problem: optimization variables' constraints value $h(x(t))$ for the different algorithms.  Exact FL-CMO refers to the unperturbed case.} 
    \label{fig:xt_perturbed_example2}
\end{figure}

\section{Conclusion}\label{sec:CON}
This paper presents a control-theoretic approach to develop new algorithms for convex and non-convex optimization with equality constraints. Based on the first-order necessary conditions, we show that the considered class of optimization problems is equivalent to a class of stabilization and output regulation problems. In this context, the Lagrange multipliers of the optimization problem serve as the control inputs for the dynamical system being regulated. We explore two methods for designing this control: Proportional-Integral (PI) control and feedback linearization.
Our analysis provides a theoretical framework for assessing the convergence properties of both methods. We rigorously prove that the PI control method converges for strongly convex problems and offer a detailed evaluation of its convergence rate. Additionally, we establish the local convergence of the feedback linearization method for non-convex problems, as well as its global convergence for strongly convex problems. Furthermore, we conduct experimental tests of both methods, validating our theoretical findings and highlighting the practical effectiveness of our proposed framework, especially in comparison to state-of-the-art optimization algorithms. 

Future work will focus on developing methods that handle inequality constraints more effectively than the approach of using squared-slack variables and non-differentiable cost functions.

\bibliographystyle{IEEEtran}
\bibliography{sophie}
\begin{IEEEbiography}[{\includegraphics[width=1in,height=1.25in,clip,keepaspectratio]{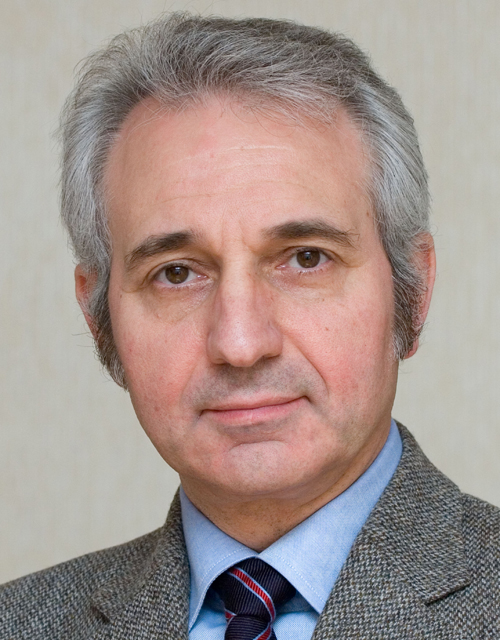}}]{Vito Cerone} received the Laurea degree in Electronic Engineering in 1984 and the PhD in System Engineering in 1989 from the Politecnico di Torino, Italy. From 1989 to 1990, he has been a Research Fellow at the School of Electronic and Electrical Engineering, University of Birmingham (UK). From 1989 to 1998, he has been an Assistant Professor with the Department of Control and Computer Engineering at the Politecnico di Torino, where he has been an Associate Professor from 1998 to 2021. Since 2021 he has been a Full Professor at the Politecnico di Torino. He is a co-author of Linear Quadratic Control: An Introduction. His main research interests include system identification, parameter estimation, optimization, and control, with applications to automotive systems.
\end{IEEEbiography}
\vspace{-0.6cm}
\begin{IEEEbiography}[{\includegraphics[width=1in,height=1.25in,clip,keepaspectratio]{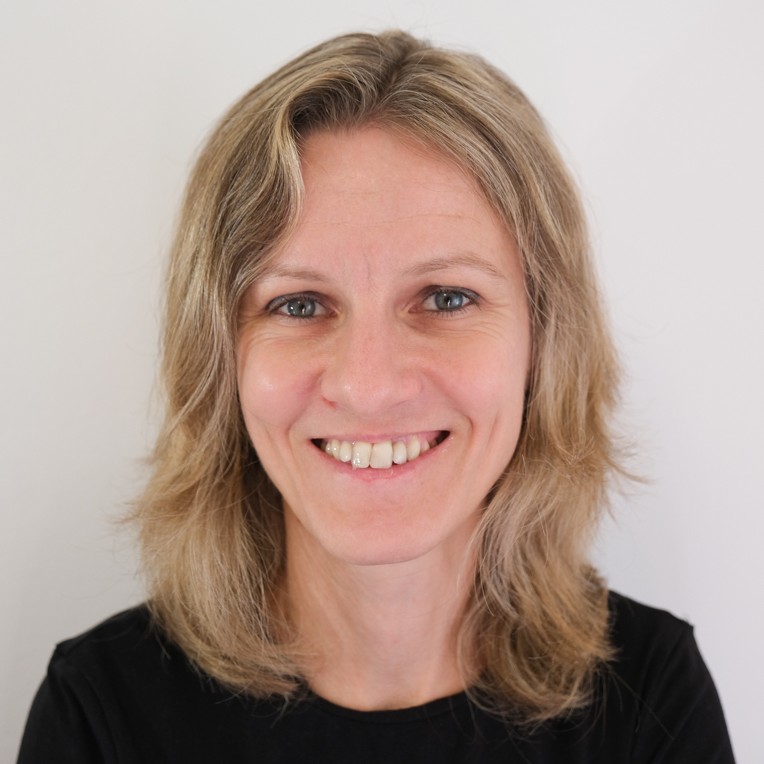}}]{Sophie M. Fosson}
(Member, IEEE) received the M.Sc. degree in applied mathematics from the Politecnico di Torino, Italy, in 2005 and
the Ph.D. degree in mathematics for the industrial technologies from Scuola Normale Superiore di Pisa, Italy, in 2011.
From 2012 to 2016, she was a Postdoctoral Associate with the Department of Electronics
and Telecommunications, Politecnico di Torino. She visited the Centre Tecnol\`{o}gic de Telecomunicacions de Catalunya, Spain, in 2013, 2014 and 2016. She was researcher at Istituto Superiore Mario Boella, Turin, Italy, in 2017. She is currently an Associate Professor with the Department of Control and Computer Engineering, Politecnico di Torino. She is Associate Editor for the IEEE Control Systems Letters. Her main research interests include sparse optimization, machine learning, system identification, control and cyber-physical systems.
\end{IEEEbiography}
\vspace{-0.6cm}
\begin{IEEEbiography}[{\includegraphics[width=1in,height=1.25in,clip,keepaspectratio]{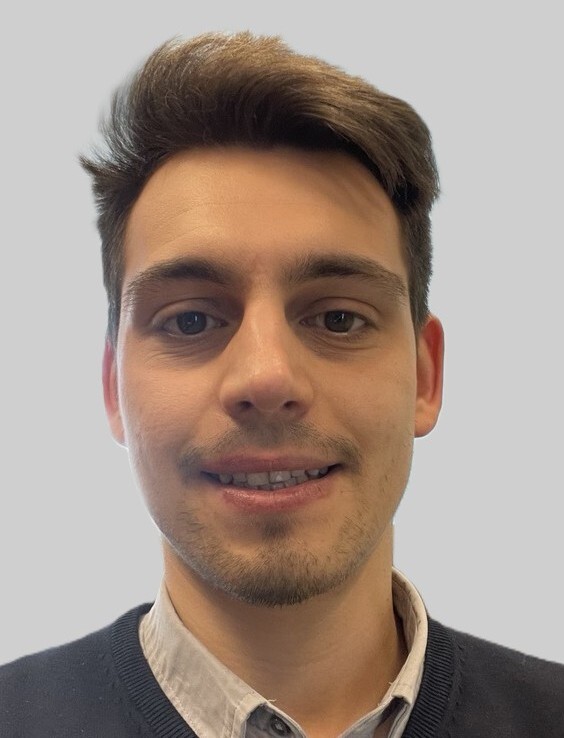}}]{Simone Pirrera } received his bachelor's degree in electronic engineering in 2019 and the master's degree in mechatronic engineering in 2021, both from Politecnico di Torino. Currently, he is a Ph.D. student in the System Identification and Control (SIC) group of the Department of Control and Computer Engineering (DAUIN) in Politecnico di Torino. His main interests are related to system identification, robust and data-driven control, mathematical optimization, and machine learning.
\end{IEEEbiography}
\vspace{-0.6cm}
\begin{IEEEbiography}[{\includegraphics[width=1in,height=1.25in,clip,keepaspectratio]{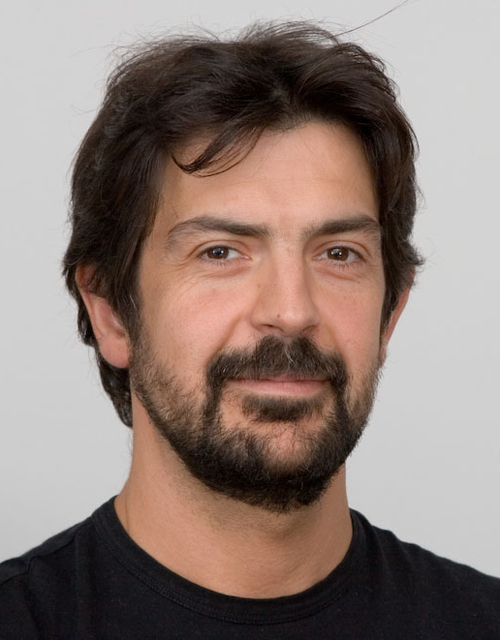}}]{Diego Regruto} received the Laurea degree in Electronic Engineering and the Ph.D. degree in System Engineering from the Politecnico di Torino, Torino, Italy, in 2000 and 2004, respectively. He is an Associate Professor of Control Theory with the Dipartimento di Automatica e Informatica, Politecnico di Torino. 

Prof. Regruto chaired the IEEE CSS Technical Committee on System Identification and Adaptive Control (TC-SIAC) from January 2013 to December 2015. He has served as Associate Editor of the IEEE Transactions on Automatic Control (from January 2016 to December 2022) and as Associate Editor of the IEEE Conference Editorial Board
(from June 2013 to December 2020). Prof. Regruto is a Senior Member of IEEE. His main research interests are in the fields of system identification, optimization and data-driven robust control methods, with application to automotive and biomedical engineering problems.
\end{IEEEbiography}

\end{document}